\documentclass{siamltex}		
\usepackage[utf8]{inputenc}
\usepackage{graphicx,multirow}  
\usepackage{verbatim,float} 
\usepackage{todonotes}
\usepackage{hyperref}
\usepackage{amssymb,amsmath}
\def\myrefs#1#2{ 
{\bigskip \noindent
{\Large \bf #2}  
 \list {[\arabic{enumi}]}{\settowidth\labelwidth{[#1]}
 \leftmargin\labelwidth 
 \advance\leftmargin\labelsep
 \usecounter{enumi} }  
 \def\newblock{\hskip .11em plus .33em minus .07em}
 \sloppy\clubpenalty4000\widowpenalty4000
 \sfcode`\.=1000\relax}  }

\def\inv{^{-1}}%
\def\nref#1{(\ref{#1})}

\def\comb#1,#2,{ \left( {#1 \atop #2 } \right)  }%
\def\prodd#1,#2,#3,{ \prod_{\scriptstyle #1 \atop\scriptstyle #2 }^{#3} }%
\def\summ#1,#2,#3,{ \sum_{\scriptstyle #1 \atop\scriptstyle #2 }^{#3} }%
\def\up#1{^{({#1})}} %

\def\CC{\mathbb{C}}

\newtheorem{algor}{{\sc Algorithm}}[section]
\newtheorem{tabl}{Table}[section]

\def\betab{\begin{tabbing} 
xxxx\=xxxx\=xxx\=xx\=xx\=xx\=xx\=xx\=xx\=xx\=xx\=xx\=xx\= \kill} 
\def\entab{\end{tabbing}\vspace{-0.12in}}

\newcommand{\eq}[1]{\begin{equation}\label{#1}}
\newcommand{\en}{\end{equation}}

\newcommand{\beeq}[1]{\begin{equation}\label{#1}}
\newcommand{\eneq}{\end{equation}}

\usepackage{graphicx,amsmath,amsfonts,amssymb,subfig}
\usepackage{color}
\usepackage[ruled,vlined,linesnumbered]{algorithm2e}

\def\cmm{{\cal M}}
\def\cma{{\cal A}}

\graphicspath{{./FIGS/}}

\title{A rational approximation method  for solving acoustic nonlinear eigenvalue problems}

\author{
  Mohamed El-Guide
  \thanks{International Water Research Institute, Mohammed VI Polytechnic University, Green City, Morocco and University of Minnesota, Department of Computer Science \& Engineering, 4-192 Keller Hall, 200 Union Street SE, Minneapolis, MN 55455, USA.
  Work supported by NSF grant 1812695.
e-mail:
\texttt{mohamed.elguide@um6p.ma}}
  \and Agnieszka Miedlar
  \thanks{University of Kansas, Department of Mathematics, 405 Snow Hall, 1460 Jayhawk Blvd. Lawrence, KS 66045-7594, USA.
    Work supported by NSF grant  1812927. e-mail:
  \texttt{amiedlar@ku.edu}} 
\and
Yousef Saad
\thanks{University of Minnesota, Department of Computer Science \& Engineering, 4-192 Keller Hall, 200 Union Street SE, Minneapolis, MN 55455, USA.
  Work supported by NSF grant 1812695.
  e-mail:
  \texttt{saad@umn.edu}}
}

\begin{document} 

\maketitle 

\begin{abstract}
  We present two approximation  methods for computing eigenfrequencies
  and   eigenmodes  of   large-scale  nonlinear   eigenvalue  problems
  resulting from   boundary element method (BEM) solutions of
   some   types   of    acoustic   eigenvalue   problems   in
  three-dimensional space.   The main idea  of the first method  is to
  approximate the  resulting boundary element matrix  within a contour
  in the  complex  plane by a high  accuracy  rational approximation  using
  the Cauchy  integral  formula.  The second  method  is  based  on
  the Chebyshev  interpolation  within  real  intervals.  A  Rayleigh-Ritz
  procedure, which  is suitable  for parallelization is  developed for
  both the Cauchy  and the Chebyshev  approximation methods when  dealing with
  large-scale  practical applications.   The  performance of  the
  proposed methods is illustrated with a variety of benchmark examples and
  large-scale  industrial applications  with  degrees of  freedom
  varying from several hundred up to around two million.
\end{abstract}

\begin{keywords} 
nonlinear eigenvalue problem, boundary element method, rational approximation, Cauchy integral formula
\end{keywords}

\section{Background and Introduction} 
The Boundary Element Method (BEM)  is a powerful approach developed to
solve integral  equations~\cite{SauS11}.  The  idea of
applying  the BEM  in many  branches  of science  and engineering  has
gained popularity over  the past few years,  e.g., in elasticity,  ground and water  flow, wave
propagation  and in  electromagnetic  problems \cite{BEM_study}.   
The most commonly  used approaches  for numerically  solving PDEs  are the
Finite Difference Method (FDM) and the Finite Element Method (FEM).  A
standard finite difference method is suitable when dealing with simple
domains (e.g. rectangular  grids), while the finite  element method can handle more complex domains. 
However, much work has to be done to numerically dicretize  a whole computational domain  (generate meshes) and  this  task   becomes  even  more  difficult   when  dealing  with
complicated domains in  higher dimensions, i.e., $d \geq  3$.  This is
where BEM becomes appealing because it allows to significatly reduce
the overall computational complexity of the solution process.  Instead
of solving a problem for the partial differential operator defined on
the  whole  domain  $\Omega$,  the boundary  element  method  uses  an
associated  boundary  integral equation  reducing  the  domain of  the
problem to the boundary $\partial \Omega$. This comes at a cost since the
matrix problem to solve in the approximation becomes dense.

In  the following,  we are  interested  in the  efficient solution  of
nonlinear  eigenvalue problems  (NLEVPs) resulting  from the  boundary
element (BE) discertization  of  the acoustic  problems.  Although  a
finite  element discretization  of  the problem  yields a  generalized
(linear) eigenvalue problem, it  requires a discretization of  the whole domain
$\Omega$  which  is  not  always  feasible, e.g.,  if  the  domain  is
unbounded.  Though the topic of NLEVPs built upon the boundary element
method (BEM) has  been around  for  a number  of years,  the lack  of
efficient  eigensolvers has  delayed a  full exploration  of BE--based
approaches. Recently,  eigenvalue solvers  based on  contour integrals
were  developed and  this made  BEM an  attractive alternative  to the
usual  contenders  when  solving  challenging   nonlinear  eigenvalue
problems~\cite{SakS03,  BEM_RRM, ZheBZZC19}.   Contour  based  methods  have  the
ability to solve NLEVPs when the  eigenvalues of interest lie inside a
given closed contour in the complex plane using rational or polynomial
approximation~\cite{EffK12, GueVBMM14, BerG17, VanBMNYBGKLNX18, LiePVM18}.  Despite these efforts, solving NLEVPs is
still  a  computationally  intensive  task.  Assembling  interpolation
matrices and  solving linear systems in the BE framework  are already
very expensive  due to the unstructured,  dense and complex nature  of the
resulting matrices.  For example, the Chebyshev  interpolation of the
BE  formulation   of  the  large-scale  accoustic   problem  discussed
in~\cite{EffK12}  results in  a generalized  eigenvalue problem  which
cannot  be easily  handled with  the state-of-the-art  linear solvers.
Another  drawback  of   this  method  is  that  the   quality  of  the
approximations   quickly  deteriorates   when  dealing   with  complex
eigenvalues.   

It  is  the  purpose  of this  paper  to  overcome  the
aforementioned  difficulties  and  develop eigensolvers  suitable  for
calculations of  eigenvalues of NLEVPs arbitrarly located in the complex plane. The  paper illustrates the performance  of the
proposed method  with a problem that arises in the  modal  analysis of
large-scale acoustic problems.

Consider the three-dimensional (3D) acoustic Helmholtz equation
\eq{eq:Helmholtz3D}
\Delta u(x)+k^2 u(x)=0,\quad x\in\Omega\subset\mathbb{R}^3,
\en
where $\Delta$ is the Laplace operator, $u(x)$ is the sound pressure at point $x$, $k=\omega/c$ is the wave number with the circular frequency $\omega$ 
and the speed of sound $c$ through the fluid medium. Equation \nref{eq:Helmholtz3D} is
subject to a homogeneous condition on its boundary $\partial \Omega$ of the form
\eq{eq:HelmBound}
a(x)u(x)+b(x)\frac{\partial u(x)}{\partial n}=0,\quad x\in\partial\Omega,
\en
where  $\frac{\partial}{\partial n}$ denotes the outward normal to the boundary at point $x$.

Using BEM yields the following Helmholtz integral equation~\cite{BEM_FEM}
\eq{eq:HelmBIM}
C(x)u(x)=\int_{\partial\Omega}\left( \frac{\partial g}{\partial n_y}u(y)-g(\|x-y\|)\frac{\partial u}{\partial n_y}\right)dy,
\en
where $C$ denotes the solid angle at point $x$, $n_y$ the surface unit normal vector at point $y$ and $g(\cdot)$ the free-space Green’s function \cite{BEM_acoustic1, BEM_acoustic2}
\eq{eq:GreenFun}
g(\|x-y\|)=\frac{e^{iz\|x-y\|}}{4\pi\|x-y\|}.
\en
The   continuous Helmholtz integral equation \nref{eq:HelmBIM}  can be discretized to  
form the following discrete problem from which the unknown boundary node values $z$ can be determined,
\eq{eq:NEVPHelmholtz}
T(z)u=0,\quad T(z):=AH(z)-BG(z),
\en
where $A$ and $B$ are diagonal matrices related to the functions $a(\cdot)$ and $b(\cdot)$ in \nref{eq:HelmBound}, $H$ and $G$ are the matrices  containing the
coefficients related to the integrals on the surface of $\frac{\partial g}{\partial n_y}$ and $g$, respectively~\cite{BEM_FEM}.
The boundary integrals are discretized by Gauss-Legendre quadrature where the
singularities of Green's function and its derivative are
isolated in the integral of revolution, and the integrations
are performed analytically using sums of elliptic integrals \cite{BEM_discretization}.
Here, $T(z) \ \in \CC^{n \times n} $ is a matrix function that is nonlinear
is $z$ and holomorphic since the free-space Green’s functions are holomorphic functions of $z$. 
Obviously, equation \eqref{eq:NEVPHelmholtz} is a NLEVP of the general form
\begin{equation}
\label{eq:NLEVPgf}
T(\lambda) u = 0,
\end{equation}
and the objective of this paper  is to develop methods for finding all
eigenvalues $z$, satisfying \eqref{eq:NEVPHelmholtz}, that are located
inside a certain  region of the complex plane enclosed by the contour
$\Gamma$.

\section{Rational and Chebyshev approximation methods for NLEVPs}
\label{sec:approx}

The first method we consider is adapted from \cite{SaaEM19} and it is based  on the  Cauchy's integral
formula. Given a Jordan curve $\Gamma$ that surrounds the eigenvalues of interest,
we express the matrix function $T(z)$ as follows:
\eq{eq:Tz}
T(z) =  \frac{1}{2 \imath \pi} \int_\Gamma \frac{T(t) }{t-z} \ dt .
\en
By replacing both occurrences of $T(\cdot)$ in \nref{eq:Tz} by
$T_{ij}(\cdot)$, one can see that the above  expression is equivalent to expressing  each individual entry
$T_{ij} (z) $ of $T(z) $ by the Cauchy integral formula. Equality \nref{eq:Tz} is valid for $z$ inside the contour
$\Gamma$ and the only requirement is that $T(z)$ be analytic inside the contour.
As is classically done \cite{GueT17} 
we use  a numerical quadrature formula to obtain the following Cauchy integral 
approximation $\widetilde T(z)$ of $T(z)$
\eq{eq:TzApp}
\widetilde T(z)\approx  \sum\limits_{i=0}^{m} \frac{\omega_i T(\sigma_i)}{z - \sigma_i} ,
\en 
where the $\sigma_i$'s are quadrature points located on the contour $\Gamma$ and the  $\omega_i$'s
the corresponding quadrature weights.

Setting $B_i=\omega_i T(\sigma_i),$ equation \nref{eq:TzApp} can be rewritten as 
\begin{align} 
  \widetilde T(z)  &= \frac{B_0}{z - \sigma_0} + \frac{B_1}{z - \sigma_1} + \ldots + \frac{B_m}{z - \sigma_m},
                     \label{eq:TzAppF0}\\
           &=   B_0 f_0(z) + B_1 f_1(z) + \ldots + B_m f_m(z) , \label{eq:TzAppF}
\end{align} 
with $f_i(z)=\frac{1}{z - \sigma_i} , i=0,\ldots,m$. 
For a given vector $u$ we now define
\[
v_i := f_i(z)u, \quad \mbox{ for } \ i=0,\ldots, m.
\]
Then the approximate nonlinear eigenvalue problem $\widetilde T(z) u = 0$ yields
\begin{equation}
\label{eq:GenNLEVP}
\widetilde T(z)u = B_0v_0 + B_1v_1 + \ldots + B_mv_m = 0. 
\end{equation}

Chebyshev interpolation of order $m$ can also be used to obtain the same form as
\nref{eq:TzAppF} of the approximation of the  matrix-valued function $T(z)$.
In this method, proposed in \cite{EffK12},
the function $T (z) $ is expanded using a degree $m$ Chebyshev polynomial expansion of the form \cite{AmiCL09}~:
\eq{eq:Chebyshev}
\widetilde T(z) = B_0\tau_0(z) + B_1\tau_1(z)+ \ldots +B_m\tau_m(z),
\en
where $B_i$ and $\tau_i(z)$ are coefficient matrices and Chebyshev basis functions, respectively. The corresponding
nonlinear eigenvalue problem is of the same form as \eqref{eq:GenNLEVP} with the vectors $v_i$ now
defined by  $v_i = \tau_i(z)u$.

The  problem \nref{eq:GenNLEVP}   for the  Cauchy interpolation, and its 
Chebyshev interpolation counterpart,  can be reformulated
as a generalized linear eigenvalue  problem: \eq{eq:sc} {\cal A} w =
\lambda {\cal M} w.  \en
For the Cauchy rational approximation we have: 
\eq{eq:sc1_cauchy}
\cma = 
\begin{bmatrix} 
  \sigma_0 I    &                &               &                 &      I          \\
                &   \sigma_1 I   &               &                 &      I           \\        
                &                & \ddots        &                & \vdots           \\
                &                &         &    \sigma_m I            &   I           \\
   -B_0         & -B_1           & \cdots        &  -B_m          &   0 
                
              \end{bmatrix}, \quad
\cmm = 
\begin{bmatrix} 
    I  &         &            &                &                   \\
       &       I &            &                &               \\        
       &         & \ddots     &                &            \\
       &         &            &   \ddots       &            \\
       &         &            &                &   0
     \end{bmatrix},
\en 
and $w = \begin{bmatrix} v_0^T, v_1^T, \ldots, v_m^T, u^T \end{bmatrix}^T$ 
whereas for the Chebyshev interpolation
\small
\eq{eq:sc1_chebyshev}
\mathcal{A} = 
\begin{bmatrix} 
	0 &I                &               &                 &               \\
	I&  0 & I  &             &                \\        
	&       \ddots        & \ddots       & \ddots                &           \\
	&                &       I  &  0            &  I           \\
	-B_0            & \cdots           & -B_{m-3}        &   C_{m-2}          &   -B_{m-1}
	
\end{bmatrix},\quad
\mathcal{M} = 2
\begin{bmatrix} 
\frac{1}{2} 	I  &     &            &                &                   \\
	&       I &          &                &               \\        
	&         & \ddots     &              &            \\
	&         &            & I        &            \\
	&         &         &                &   B_m
\end{bmatrix},
\en
\normalsize 
where $C_{m-2}  \equiv B_m-B_{m-2}$ 
and $w = \begin{bmatrix} u^T, v_1^T, \ldots, v_{m-1}^T, v_m^T\end{bmatrix}^T$.

With regards  to the rational  approximation described above,  we note
that  an alternative  that  has been  used with  some  success in  the
literature     is     the    Barycentric     approximation     formula
\cite{BerG17}. However, our tests with this technique showed  no significant  improvement
in our context  over  the simple  Cauchy
formula used above. Note that it is also possible to exploit other polynomials,
using different classes of orthogonal polynomials but we will restrict our
attention to Chebyshev polynomials of the first kind. 
Finally note that Chebyshev approximation works best for eigenvalues located in an interval while
the Cauchy rational approximation is suitable for general complex spectra.

\section{Rayleigh–Ritz procedure for BEM eigenvalue problem}\label{RR}
Let $U$ be a basis of dimension $\nu$ of a subspace that contains good approximations of the eigenvectors of
the NLEVP problem \nref{eq:NEVPHelmholtz}.
Then, it is possible to  apply a Rayleigh-Ritz procedure to \nref{eq:NEVPHelmholtz} to obtain
approximate eigenpairs. The  approximate eigenvector will be of the form
$u=U y $ \ with \ $ y \in \mathbb{C}^{\nu}$. Then expressing that $T(z) u $ is orthogonal to the range of
$U$ yields the projected problem $U^H T_U(z)u = 0$ or,
\eq{eq:RR_approx}
 \widehat{B}_0f_0(z) y +\widehat{B}_1f_1(z)  y  + \ldots + \widehat{B}_mf_m(z)  y  =0,
\en
where $\widehat{B}_i=U^H B_iU$.
We will denote $T_U(z) $ the projected operator, namely, 
\eq{eq:reduced_approx}
T_U(z)= \widehat{B}_0f_0(z)+\widehat{B}_1f_1(z)+ \ldots +\widehat{B}_mf_m(z). 
\en

Then,  applying the same procedure as before to the projected problem  we see that
\nref{eq:RR_approx} becomes: 
\eq{eq:RR_Cauchy}
\widehat{B}_0\widehat v_0+\widehat{B}_1\widehat v_1+\ldots+\widehat{B}_m \widehat v_m=0,
\en
with $\widehat v_i = \frac{ y }{z - \sigma_i}$ in the case of rational approximation and $\widehat v_i = \tau_i(z) y $
when a Rayleigh-Ritz procedure is applied to the Chebyshev interpolation.

\subsection{Solution of the reduced NLEVP}
Analogously to what was discussed in  Section~\ref{sec:approx}, the problem  \nref{eq:RR_Cauchy} for the
Cauchy rational approximation, as  well as its Chebyshev interpolation
counterpart, can  be written down  in a block form  \eqref{eq:sc}, but
now of much smaller dimension.
The projected nonlinear problem \eqref{eq:RR_Cauchy} yields the following linearized problem
\begin{equation}
\label{eq:EVP_red}
\cma w = \lambda \cmm w,
\end{equation}
with $w = \begin{bmatrix} \widehat v_0^T, \widehat v_2^T, \ldots, \widehat v_m^T,  y ^T\end{bmatrix}^T$ and
\eq{eq:sc1_cauchy_red}
\mathcal{A} = 
\begin{bmatrix} 
  \sigma_1 I    &                &               &                 &      I          \\
                &   \sigma_2 I   &               &                 &      I           \\        
                &                & \ddots        &                & \vdots           \\
                &                &         &    \sigma_m I            &   I           \\
   -\widehat B_1         & -\widehat B_2           & \cdots        &  -\widehat B_m          &   0 
                
              \end{bmatrix}, \quad
\mathcal{M} = 
\begin{bmatrix} 
    I  &         &            &                &                   \\
       &       I &            &                &               \\        
       &         & \ddots     &                &            \\
       &         &            &   \ddots       &            \\
       &         &            &                &   0
     \end{bmatrix}
\en
for the Cauchy rational approximation, and
$w = \begin{bmatrix}  y ^T, \widehat v_1^T, \ldots, \widehat v_{m-1}^T, \widehat v_m^T\end{bmatrix}^T$,
\small
\eq{eq:sc1_chebyshev_red}
\mathcal{A} = 
\begin{bmatrix} 
	0 &I                &               &                 &               \\
	I&  0 & I  &             &                \\        
	&       \ddots        & \ddots       & \ddots                &           \\
	&                &       I  &  0            &  I           \\
	-\widehat{B}_0            & \cdots           & -\widehat{B}_{m-3}        & \widehat{C}_{m-2}
                 &   -\widehat{B}_{m-1}
\end{bmatrix}, \ \                
\mathcal{M} = 
\begin{bmatrix} 
	I  &     &            &                &                   \\
	&       2I &          &                &               \\        
	&         & \ddots     &              &            \\
	&         &            & 2I        &            \\
	&         &         &                &   2\widehat{B}_m
\end{bmatrix}
\en
\normalsize
for      the     Chebyshev      interpolation     where
$\widehat{C}_{m-2}  = \widehat{B}_m-\widehat{B}_{m-2}$.   If $\nu$  is
fairly small,  the problem  \nref{eq:EVP_red} can be  solved directly, i.e.,
using standard dense packages. 
When $\nu  $ is larger,  the problem must  be handled differently by some iterative procedure. 
Since for  BEM  problems the  matrices  $B_i$ are  generally
complex,  dense  and  unstructured, solving  these  linear  eigenvalue
problems  can be  computationally  expensive.  Therefore,  it may be
advantageous to rely on subspace iteration or an Arnoldi-type method to solve
\eqref{eq:EVP_red}.

Note that for both the Cauchy rational approximation and the Chebyshev
interpolation method, the matrices $\cma$  and $\cmm$ don't have to be
formed   explicitly.  If   the   partial  solution   of  the   problem
\eqref{eq:EVP_red}  are  of interest,  effective methods  such as  the
Implicitly Restarted Arnoldi method can  be used to find a few of the extremal
eigenvalues.  Unfortunately,  these   methods  become  expensive  when
the eigenvalues of interest are deep inside the spectrum.

Alternatively, we can  solve the interior eigenvalue  problem with the
help of the shift-and-invert technique, which replaces the solution of
the generalized eigenvalue problem \eqref{eq:EVP_red} by the following
problem \eq{eq:sc2} \mathcal{ H} w = \frac{1}{\lambda-\sigma} w,\qquad
\mathcal{ H} := \left({\cal A  }-\sigma {\cal M} \right)^{-1}{\cal M}.
\en
Using the  Arnoldi or the  subspace iteration  method to  extract extremal
eigenvalues  of \nref{eq:sc2}  will  result in  approximations of  the
eigenvalues  of \eqref{eq:EVP_red}  closest to  $\sigma$.  Again,  the
matrix $\mathcal{H}$ need not be formed explicitly to compute
the  matrix-vector product  $y=\mathcal{H}x$.  Instead, we  can use  a
simple  $LU$ factorization  that takes  advantage of  the sparsity  of
${\cal  A   }$  and  ${\cal  M   }$.  First,  note  that   the  matrix
$\left({\cal A }-\sigma {\cal M} \right)$ is of the form \eq{eq:sc3}
 \begin{bmatrix} 
 	D & F \\
 	B & C
 \end{bmatrix}.
 \en
 By exploiting  the sparsity of the matrices $D$ and $F$, we can easily form the following $LU$ factorization
 \eq{eq:sc4}
 L=\begin{bmatrix} 
 	I & 0 \\
 	BD^{-1} & I
 \end{bmatrix}, \qquad 
 U=\begin{bmatrix} 
 	D & F \\
 	0 & S
 \end{bmatrix},
 \en
 where $S=C-BD^{-1}F$ is known  as the \emph{Schur complement} of the block $C$. 
 With the use of matrix $S$, we can use the Arnoldi algorithm on vectors of shorter length.
 Solving the shifted and inverted problem \nref{eq:sc2} with Arnoldi algorithm
 requires solving linear systems of the form
 \eq{eq:sc5}
 \begin{bmatrix} 
 	D & F \\
 	B & C
 \end{bmatrix}\begin{bmatrix} 
 	x\\
 	y
 \end{bmatrix}=\begin{bmatrix} 
 	a\\
 	b
 \end{bmatrix}.
 \en
 Using the Schur complement $S$, $y$ can be easily obtained by solving $Sy=b-BD^{-1}a$ and since $D$ is a diagonal matrix and $F$ is a block of identity matrices, one can determine $x$ by using  the relation $Dx+Fy=a$.

\subsection{Construction of the subspace of approximants}
We begin this section by noting that the Arnoldi-type  or subspace  iteration  methods
discussed in  the previous  section can   be  applied to  a linear
eigenvalue problem  $\cma w =  \lambda \cmm w$ obtained  directly from
\eqref{eq:GenNLEVP}.  However, proceeding in this way would 
require either solving linear eigenvalue  problems of size $mn+n$ when
using Arnoldi-type methods, or storing vectors of length $mn+n$ in the
subspace  iteration method, and this can be  computationally   expensive
when $m$ is large. Therefore,  it is important to develop
a technique  that allows  to work with subspaces of  smaller dimensions that
requires storing shorter  vectors. A procedure of this type, which works with
subspaces of dimention $m$ is presented next.

Let  us first  consider  a  large linear  eigenproblem  of the  form
\nref{eq:sc} obtained  from \eqref{eq:GenNLEVP} without  a projection.
To introduce the approach that works with vectors of dimension $n$, we first point
out that for an approximate eigenpair $(\lambda,u)$, $u$ is the bottom
(resp. top part) of an approximate eigenvector $w$ of the large linear
eigenvalue problem \eqref{eq:sc} associated with \eqref{eq:TzAppF} for
the Cauchy rational approximation, (resp. \eqref{eq:Chebyshev} for the
Chebyshev interpolation).  Let $W^{(0)}$  be a random initial set of
$\nu$ basis vectors of a certain subspace,
where each of the $\nu$ columns  is of the form \footnote{Here we use
  Matlab notation:  $[v; u]$  is a  vector that stacks  $v$ on  top of
  $u$.}  $w=[v; u]$  (resp.  $w=[u; v]$) for  the splitting associated
with the  Cauchy rational approximation (resp. Chebyshev interpolation).   Next, in
order to make  these initial random vectors close  to the eigenvectors
of  interest, we  apply $q$  steps of  the inverse  power method  with
matrix  $\mathcal{H}$ in  \nref{eq:sc2}  to each  column of  $W^{(0)}$
separately.   A  subspace  of  dimension  $n$  that  approximates  the
eigenvectors of  \nref{eq:GenNLEVP} is  then obtained from  the bottom
parts  (resp.  top  parts) of  the  processed columns  for the  Cauchy
rational approximation (rep.  Chebyshev interpolation).  Although this
process  involves the  column vectors  of $W^{(0)}$,  only vectors  of
length $n$  need to be  saved and the  iterates $v$ can  be discarded.
The  accuracy of  the extracted  eigenpairs obtained  from applying  a
Rayleigh-Ritz projection can  be further refined by updating  $U$ in a
process  that takes  advantage  of the  structure  of the  approximate
eigenvectors.   Let  $(\lambda,u)$  be  an  approximate  eigenpair  of
\nref{eq:GenNLEVP} obtained  from applying a  Rayleigh-Ritz projection
using $U$.  The new redefined vector $w$ for each interpolation method
is  discussed  next.  For  the  Cauchy  rational approximation, the  vector
$v  = [v_1; v_2; ...; v_m]^T$  (the  top part  of  vector $w$),  is
obtained by  setting $v_i=\frac{u}{\lambda  - \sigma_i}$,  whereas for
the Chebyshev interpolation (the bottom part of vector $w$) it is defined
by setting $v_i=\tau_i(\lambda)u$.
 
\subsection{The inverse power method}
The straightforward linearizations \eqref{eq:sc1_cauchy} of the Cauchy
rational approximation  and \eqref{eq:sc1_chebyshev} of  the Chebyshev
interpolation,    discussed   in    Section   \ref{sec:approx},    are
high dimensional problems and they become computationally demanding 
as  the order $m$ of  the approximations  grows.  The  Rayleigh-Ritz
approach discussed above is inexpensive even if $m$ is large. The biggest
computational task  of the presented Rayleigh-Ritz
projection lies in performing $q$  steps of the  inverse power  method with
the matrix $\left({\cal A  }-\sigma {\cal M} \right)^{-1}{\cal  M}$. It is
the purpose of  the following discussion to show how  each step of the
inverse  power  method  can  be carried  out  inexpensively.  For 
simplicity, we  will assume that the  shift $\sigma$ is the  center of
the  unit circle  (resp. interval  $[-1,1]$) for  the Cauchy  rational
approximation  (resp. Chebyshev  interpolation).   This  is a  natural
choice, since  any circle in  the complex plane  can be scaled  to the
unit  circle and  any  real  interval $[a,b]$  can  be  scaled to  the
interval $[-1,1]$.  Throughout this  discussion, the superscript $j$ will
correspond  to the  iteration  number, while  the  subscript $i$  will
correspond  to  the blocks  of  the  vectors  $v^{(j)}$. We  begin  by
discussing  the   inverse  power   method  for  the   Cauchy  rational
approximation.


\paragraph{Inverse power iteration for the Cauchy rational approximation}
 
%
%
For the Cauchy interpolation, each step of the inverse power iteration method requires solving a linear system
\eq{eq:lineaCauchy}
\mathcal{A}w\up{j+1}= y\up{j} \ \ \mbox{ with } \ \ y\up{j}=\mathcal{M}w\up{j} \ \mbox{ and } \ w\up{j}=[v\up{j}; u\up{j}], 
\en 
which is of the form \nref{eq:sc5}. Therefore, the iterates of the inverse power method can be determined by solving 
 \begin{align}
 S u\up{j+1} & =b, \ \mbox{ with } \ b=\big(u\up{j} - B D\inv v\up{j}\big),  \label{eq:wjp1}\\
 D v\up{j+1} & =(v\up{j} - F u\up{j+1}) .  \label{eq:vjp1}
 \end{align}
 Since $D$ is a diagonal matrix and $F$ is a block vector of identity matrices, $v_i\up{j+1}$ are determined by
 \eq{eq:viupd}
 v_i\up{j+1} = \frac{v_i \up{j} - u\up{j+1}}{\sigma_i},\quad  i=0,\ldots, m.
 \en
 Again, exploiting the structure of $D$ and $F$, the iterate $u\up{j+1}$ can be obtained by solving (\ref{eq:wjp1}) with
 \eq{eq:linsyssim}
 S = -\sum_{i=0}^m\frac{B_i}{\sigma_i}, \quad \text{and} \quad b=u\up{j}-\sum_{i=0}^m\frac{B_i}{\sigma_i}v_i\up{j}.
 \en
Algorithm~\ref{alg:InvItCauchy} performs one step of the inverse power iteration for the Cauchy rational approximation. 
 
\IncMargin{1em}
\begin{algorithm}[h!]
	\SetKwFunction{mgs}{mgs}
	\SetKwInOut{Input}{Input}
	\SetKwInOut{Output}{Output}
	\Input{$D, F, B$ and $C = 0$ as defined in \eqref{eq:sc5}, $w^{(j)} = \begin{bmatrix} v^{(j)} \\ u^{(j)} \end{bmatrix}$}
	\Output{$w^{(j+1)} = \begin{bmatrix} v^{(j+1)} \\ u^{(j+1)} \end{bmatrix}$}
	\BlankLine
	Compute $b = u^{(j)} - BD^{-1}v^{(j)} = u^{(j)} - \sum\limits_{i=0}^{m}\frac{B_i}{\sigma_i}v_i^{(j)}$\;
	Solve $Su^{(j+1)} = b$, with the Schur complement matrix $S = C - \sum\limits_{i=0}^{m} \frac{B_i}{\sigma_i} = - \sum\limits_{i=0}^{m} \frac{B_i}{\sigma_i}$\;
	Set $v_i^{(j+1)} = \frac{v_i^{(j)} - u^{(j+1)}}{\sigma_i}$\;
	\Return{$v^{(j+1)}$, $u^{(j+1)}$}
	\caption{One step of inverse power method for Cauchy approximation}
	\label{alg:InvItCauchy}
\end{algorithm}\DecMargin{1em}

\paragraph{Inverse power iteration for the Chebyshev interpolation}
Recall that the iterates obtained from the inverse power method
for    the    Chebyshev    interpolation    can    be    written    as
$w\up{j}  = [u\up{j};  \ v\up{j}  ] $.   Similarly to the  Cauchy
rational approximation, each step of the inverse power method requires
solving        the        linear       system
\eq{eq:lineaCheb}
\mathcal{A}w\up{j+1}=y\up{j},    \    \    \mbox{   with    }    \    \
y\up{j}=\mathcal{M}w\up{j}.
\en
Since     $\mathcal{M}$     is     a    block     diagonal     matrix,
$y\up{j}=\mathcal{M}w\up{j}$  can  be  easily evaluated.
The  question that  remains to  be answered  is how  to solve
efficiently  the  linear  system  $\mathcal{A}w\up{j+1}=y\up{j}$.   By
taking  advantage of  the  block structure  of  $\mathcal{A}$ for  the
Chebyshev interpolation, it follows naturally that this problem can be
treated    by   performing    the   following    steps,
see~\cite[Section 2.3]{EffK12}.
To compute the bottom part $v\up{j+1}$ of $w\up{j+1}$ we
will use the recursion
\eq{eq:reccurenceCheb_odd}
 v_1\up{j+1}=y_0\up{j},\quad v_{2i+1}\up{j+1}=y_{2i}\up{j}-v_{2i-1}\up{j+1}, \quad i=1,2,\ldots,
\en 
for  odd-numbered blocks and
\eq{eq:temp1}
v_0^{(j+1)} = u^{(j)} , \quad v_{2i}^{(j+1)} = y_{2i}^{(j)} - v_{2i-2}^{(j+1)}, \quad i = 1,2, \ldots.
\en
%
for even-numbered blocks.
Since the blocks $v_{2i-2}^{(j+1)}$ in \eqref{eq:temp1} are even-numbered, we can further expand the recurrence relation, i.e.,
\eq{eq:reccurenceCheb_even}
v_0^{(j+1)} = u^{(j)}, \quad v_{2i}^{(j+1)} = \widehat y_{2i-1}^{(j+1)} + (-1)^{i}v_0^{(j)}, \quad i = 1,2, \ldots,
\en
where
\[
\widehat y_1^{(j+1)} = y_1^{(j)}, \quad \widehat y_{2i+1}^{(j+1)} = y_{2i+1}^{(j)} - \widehat y_{2i-1}^{(j+1)}, \quad i=1,2,\ldots \ .
\]

Since $v_0 = \tau_0(z)u$ and $\tau_0 = 1$ (zeroth Chebyshev polynomial), $u\up{j+1}$ is the top part of vector $w^{(j+1)}$, i.e., $u^{(j+1)} = v_0^{(j+1)}$
and it can be obtained by solving 
\eq{eq:RedCheb}
 Gu\up{j+1}=b.
\en 
Given the number of quadrature nodes $m$, let us consider the Euclidean division of $m$ by $2$, i.e.,
$m = 2\cdot q + r$, with quotient $q$ and remainder $r$. Then the matrix $G$ has the following form
\begin{equation}
 \label{eq:Gmtx}
G = \sum_{i=0}^{q}(-1)^{i+1}B_{2i}. 
\end{equation}
The  vector $b$ depends on the parity of $m$.
If $m$ is odd
\begin{equation}
\label{eq:bodd}
b = \sum\limits_{i=0}^{q-1}B_{2i+1}v_{2i+1}\up{j+1} + \sum\limits_{i=1}^{q}B_{2i}\widehat y_{2i-1}\up{j+1} + y_{m-1}\up{j} - B_m \big( \sum\limits_{i=0}^{q-1}(-1)^{q-i} y_{2i}\up{j} \big),
\end{equation}
and when it is even, then 
\begin{equation}
\label{eq:beven}
b = \sum\limits_{i=0}^{q-1}B_{2i+1}v_{2i+1}\up{j+1} + \sum\limits_{i=1}^{q-1}B_{2i}\widehat y_{2i-1}\up{j+1} + y_{m-1}\up{j} - B_m \big( \sum\limits_{i=0}^{q-1}(-1)^{q-i}y_{2i+1}\up{j} \big) .
\end{equation}

Algorithm~\ref{alg:InvItChebyshev} implements  one step of inverse power method for the Chebyshev interpolation.
%
%

\IncMargin{1em}
\begin{algorithm}[h!]
	\SetKwFunction{mgs}{mgs}
	\SetKwInOut{Input}{Input}
	\SetKwInOut{Output}{Output}
	\Input{$B_0, \ldots, B_m$, $w^{(j)} = \begin{bmatrix} v^{(j)} \\ u^{(j)} \end{bmatrix}$}
	\Output{$w^{(j+1)} = \begin{bmatrix} v^{(j+1)} \\ u^{(j+1)} \end{bmatrix}$}
	\BlankLine
	Compute $v^{(j+1)}$ using recurences \eqref{eq:reccurenceCheb_odd} and \eqref{eq:reccurenceCheb_even}\;
	Form matrix $G$ defined in \eqref{eq:Gmtx} and right-hand side vector $b$ using \eqref{eq:bodd} or \eqref{eq:beven}\;
	Solve linear system $Gu^{(j+1)} = b$\;
	\Return{$v^{(j+1)}$, $u^{(j+1)}$}
	\caption{One step of inverse power method for Chebyshev approximation}
	\label{alg:InvItChebyshev}
\end{algorithm}\DecMargin{1em}

To this end, only one LU factorization is required -- of the Schur complement matrix $S$ in the case of the Cauchy rational approximation
or matrix $G$ for the Chebyshev interpolation -- in the preprocessing step for all $q$ steps of the inverse power method.

\IncMargin{1em}
\begin{algorithm}[h!]
	\SetKwFunction{mgs}{mgs}
	\SetKwInOut{Input}{Input}\SetKwInOut{Output}{Output}
	\Input{ 
          Subspace dimension $\nu$; $q$; 
          Number of eigenvalues $k$ (with  $k\le \nu$) } 
        \Output{$\lambda_1,   \ldots,  \lambda_k$, $U_k$} 
        \For{$j=1:\nu$}{
           Select $w = [v; \ u]$ (or $w =[u; \ v])$ \tcc*[r]{Initially random vectors} 
           Run $q$ steps of Algorithm~\ref{alg:InvItCauchy} or \ref{alg:InvItChebyshev} starting  with $w$
           \;
           If $w = [v;\  u]$ (or $w = [u; \ v]$) is the last iterate, then set $U(:,j)  = u$\;
         }
         Use  $U$   to  compute $\widehat B_i, \ 0=1,...,m$  from
        \eqref{eq:RR_Cauchy}\;
        Solve   the   reduced   eigenvalue   problem
        \nref{eq:EVP_red} associated with \nref{eq:sc1_cauchy_red} or \eqref{eq:sc1_chebyshev_red}\;
	\Return{$\lambda_1, \ldots, \lambda_k$ and eigenvector matrix $U_k$}
	\caption{Reduced subspace iteration (no restarts) for Cauchy (or Chebyshev) approximation} 
	\label{alg:subsit}
\end{algorithm}\DecMargin{1em}



\section{Numerical Experiments} \label{sec:exper} 
This section will illustrate the behavior of the  approaches presented in this paper for solving
nonlinear eigenvalue problems in the form \eqref{eq:NEVPHelmholtz}
resulting from boundary element disretization of \eqref{eq:Helmholtz3D} -- \eqref{eq:HelmBound}.
All   experiments  were performed with {\sc Matlab} R2018a. Furthermore, 
computations in Example 3 were performed in parallel on a Linux cluster
at the Minnesota Supercomputer Institute that has  $32$ cores and $31.180$ GB per-core memory.

For the presented examples, the contour $\Gamma$
is either circular or elliptic and the eigenvalues of interest are those
closest to the center of $\Gamma$, i.e.,
in Algorithms \ref{alg:InvItCauchy} and \ref{alg:InvItChebyshev} the shift $\sigma$ is selected to be the center of the
region enclosed by the contour $\Gamma$.

In the case of a circular contour,  the $m$ quadrature nodes and weights
used  to  perform  the  numerical integration  to approximate  the
functions $f_j$   inside the  contour $\Gamma$ were  generated using 
the Gauss-Legendre quadrature rule. To illustrate the
effectiveness of  the proposed approaches, we  compare the eigenvalues
obtained by each algorithm either with exact eigenvalues or the approximations obtained 
by the Beyn's method~\cite{Bey12} or/and via a corresponding linearization.

\subsection*{Example 1}
As our first example, we consider the 3D Laplace eigenvalue problem \nref{eq:Helmholtz3D} on the unit
cube $\Omega = [0, 1]^3$ with homogeneous Dirichlet boundary conditions, i.e., \eqref{eq:HelmBound} with $a(x)=1$ and $b(x)=0$.
The exact eigenvalues for this problem are known and given by
\eq{eq:cubeev}
k=\sqrt{n_1^2+n_2^2+n_3^2}, \quad  n_i=0,1,2, \ldots \ .
\en
We are interested in the six smallest eigenvalues (including multiplicities) of \nref{eq:Helmholtz3D} presented in Table~\ref{tab:LaplaceEvals}.

\begin{figure}
	\centering
	\includegraphics[width=0.49\textwidth]{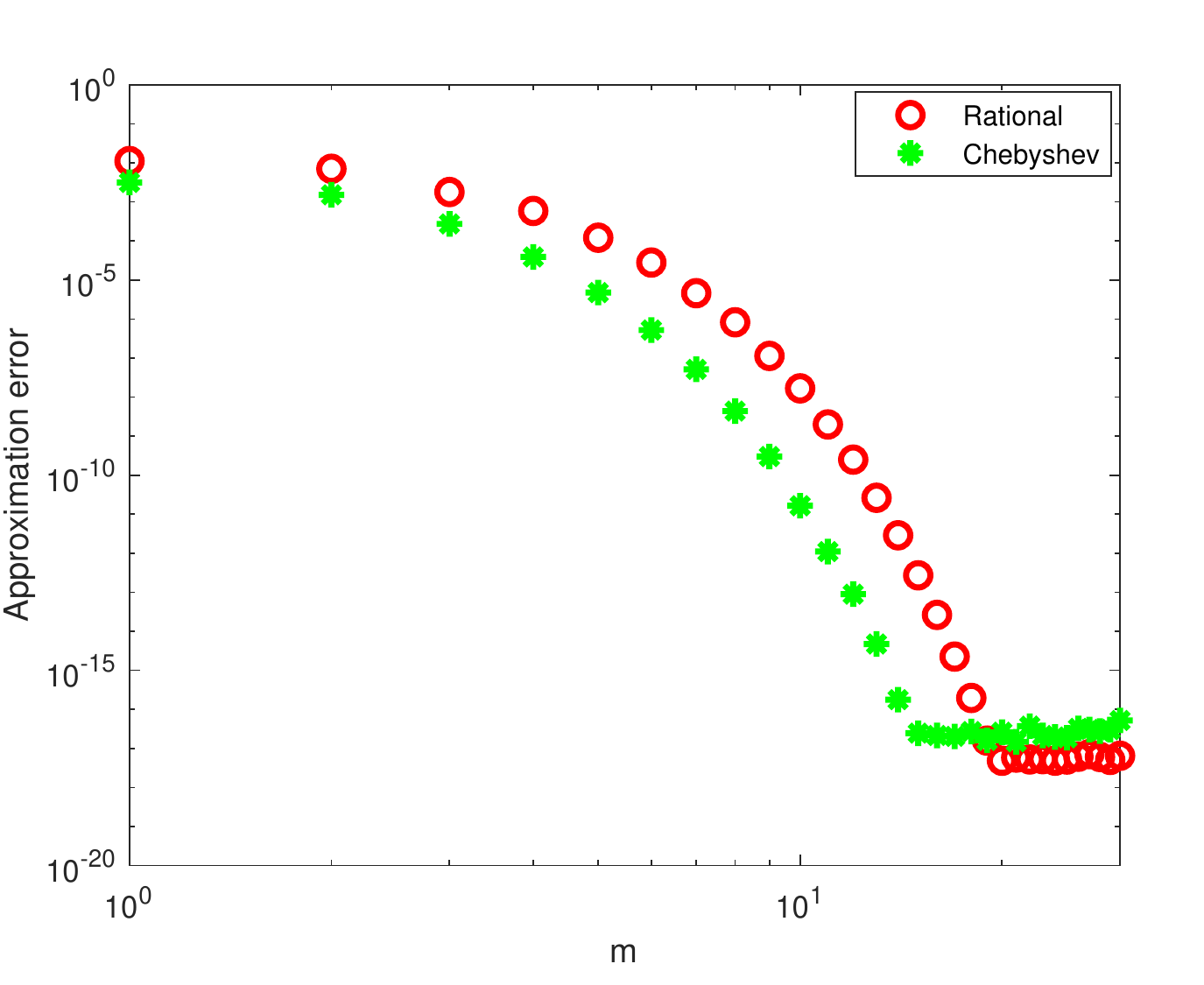}
	\includegraphics[width=0.49\textwidth]{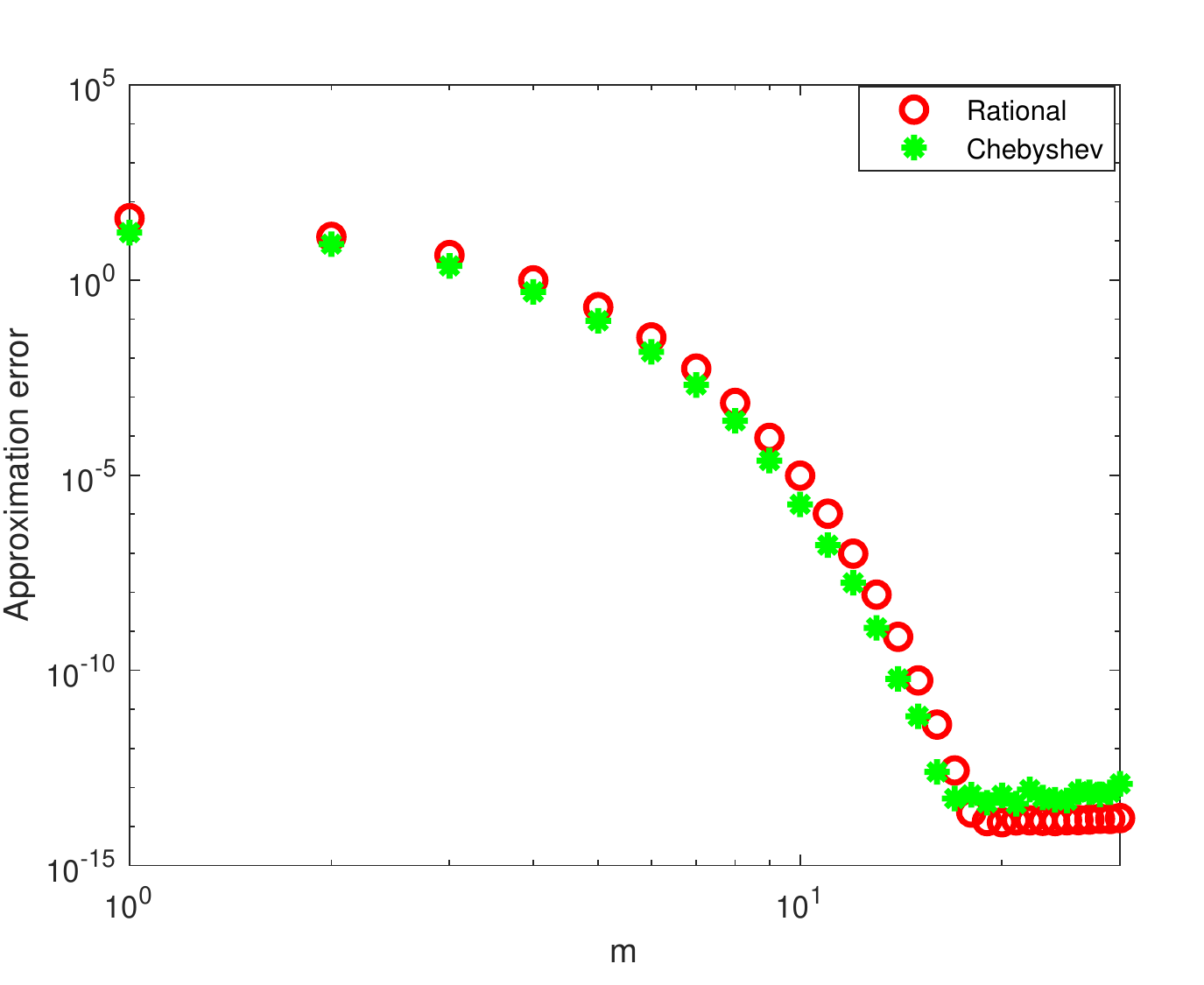}
	\caption{Left:  Approximation error versus the order of the approximation $m$ inside a unit circle. 
	         Right:  Approximation error versus the order of the approximation $m$ inside an ellipse centered at $c = 0$ with semi-major axis
		$r_x = 1$ and semi-minor axis $r_y = 0.2$.}\label{fig:RatChebApp}
\end{figure}

\begin{figure}
	\centering
	\includegraphics[width=0.49\textwidth]{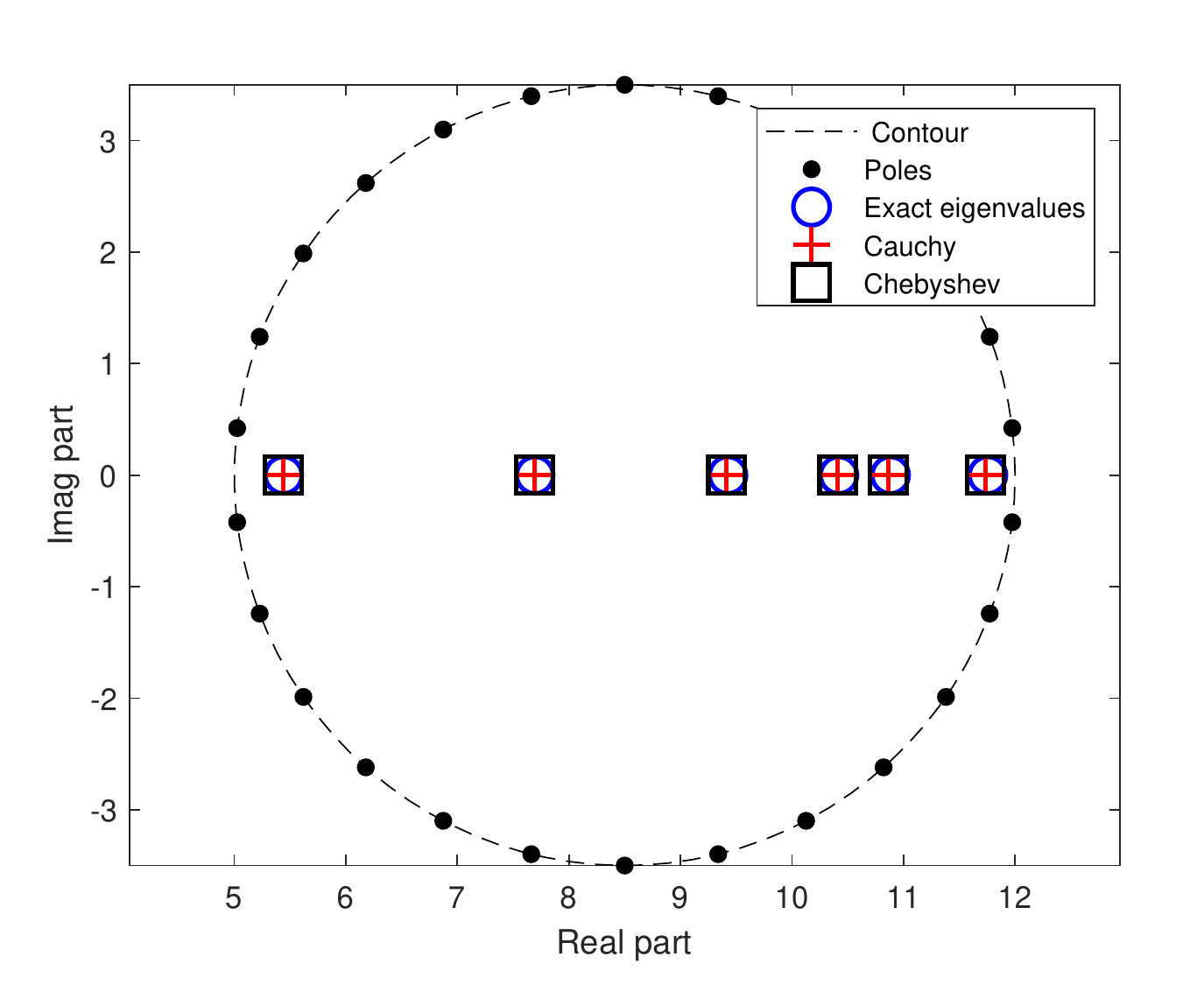}
	\includegraphics[width=0.49\textwidth]{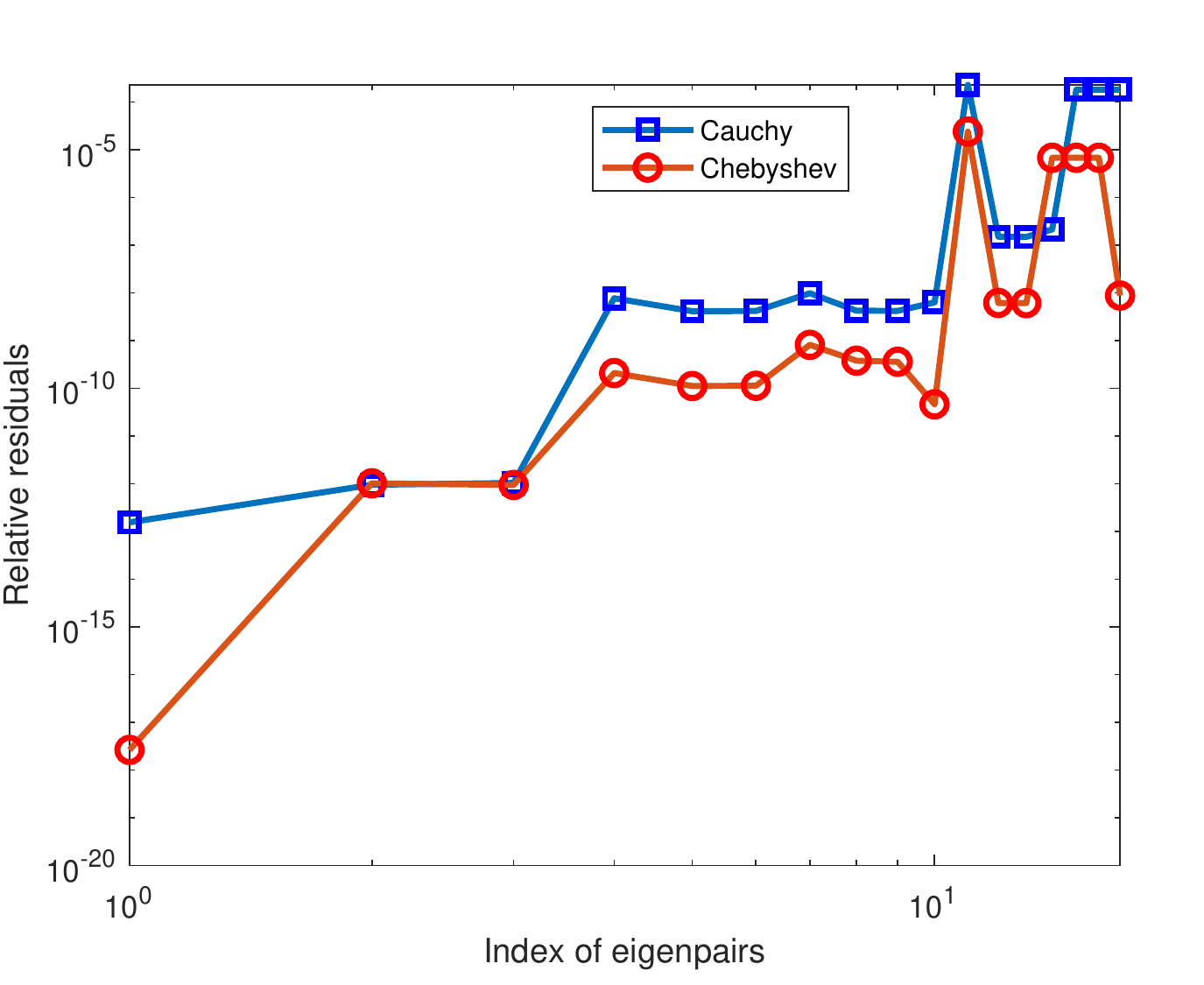}
	\caption{Left: The eigenvalues of \nref{eq:Helmholtz3D} with homogeneous Dirichlet boundary conditions inside a
		circle centered at $c = 8.5$ with radius $r = 3.5$ (circles) computed via \nref{eq:sc} (plus) and
		Chebyshev interpolation method inside the real interval $[5,12]$ (squares).
		Right:  The relative residuals $\|T(\lambda)u\|_2/\|u\|_2$ of the computed eigenpairs.}\label{fig:Cube_Dirichlet}
\end{figure}

\begin{figure}
	\centering
	\includegraphics[width=0.49\textwidth]{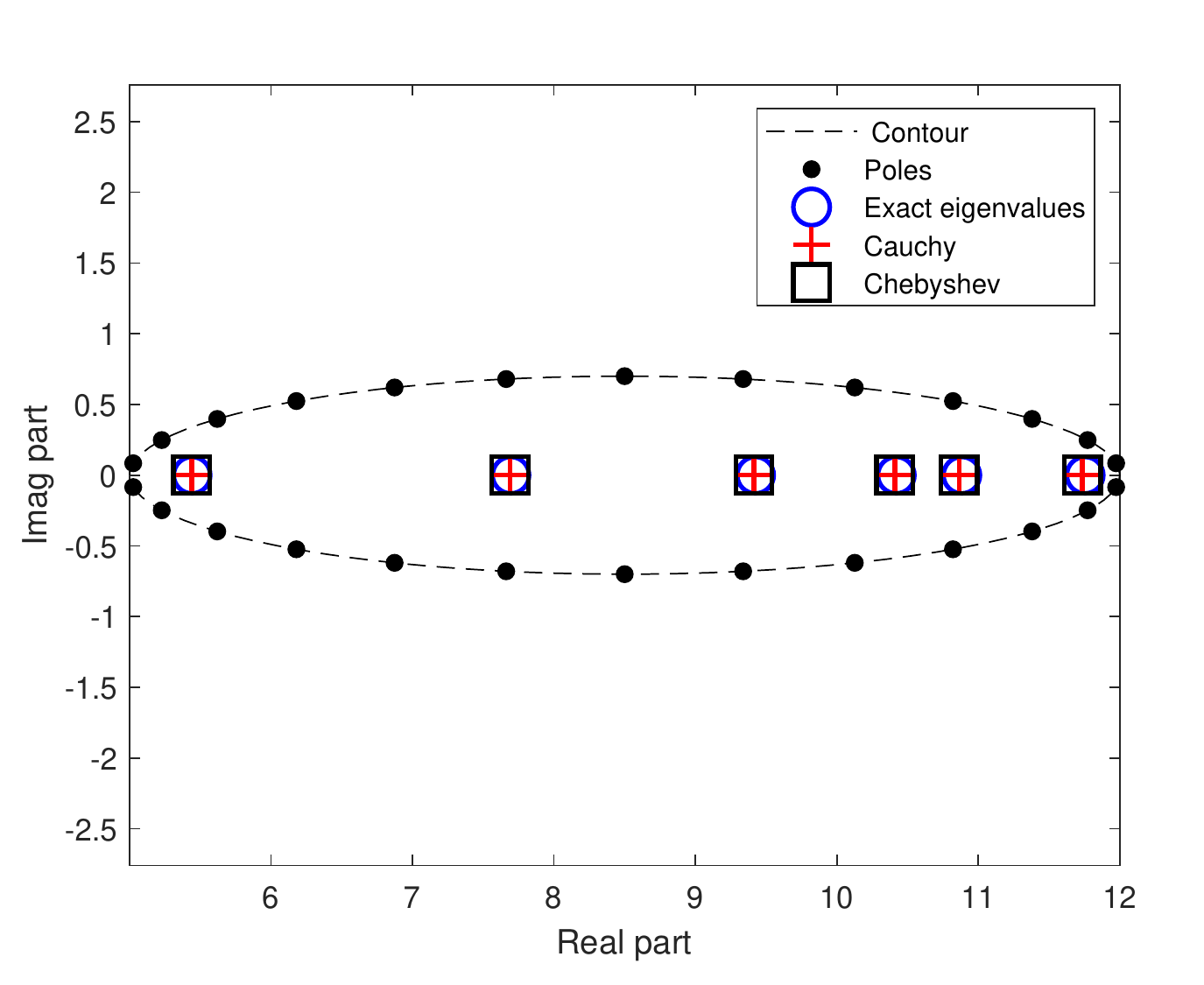}
	\includegraphics[width=0.49\textwidth]{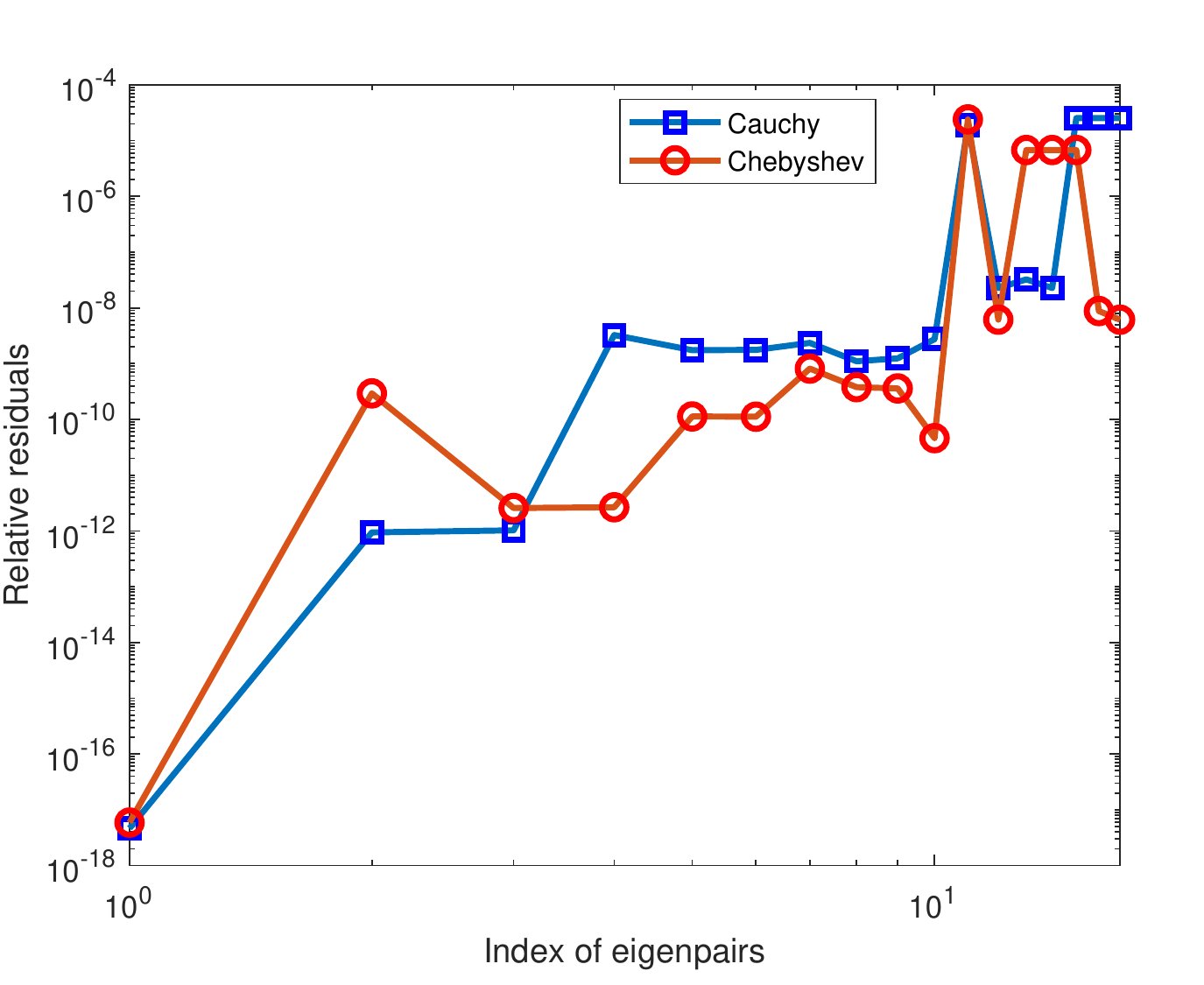}
	\caption{Left: The eigenvalues of \nref{eq:Helmholtz3D} with homogeneous Dirichlet boundary conditions inside an ellipse centered
		at $c = 8.5$ with semi-major axis $r_x = 4.5$ and semi-minor axis $r_y = 0.2$ (circles) computed 
	        computed via \nref{eq:sc} (plus) and Chebyshev interpolation method inside the real interval $[5,12]$ (squares).
	        Right:  The relative residuals $\|T(\lambda)u\|_2/\|u\|_2$ of the computed eigenpairs.}\label{fig:ellipse_BEM}
\end{figure}

\begin{table}[H]
\begin{center}
\begin{tabular}{c c c}
\hline
 no. & eigenvalue & multiplicity\\
\hline
1 & 5.441398  & 1 \\
2 & 7.695299  & 3 \\
3 & 9.424778  & 3 \\
4 & 10.419484 & 3 \\
5 & 10.882796 & 1 \\
6 & 11.754763 & 6 \\
\hline\\
\end{tabular}
\end{center}
\caption{Approximations of the $6$ smallest eigenvalues (including multiplicities) of the 3D Laplace eigenvalue problem on $\Omega = [0,1]^3$ with homogeneous Dirichlet boundary conditions~\cite[Table 1]{EffK12}.}
\label{tab:LaplaceEvals}
\end{table}
To  determine   these  eigenvalues  using  the   Cauchy  approximation
technique,   we  will   build  the   rational  approximation   of  the
matrix-valued   function  $T(\cdot)$   using  circular   and  elliptic
contours. First, we  compare the accuracy between  the Cauchy rational
approximation  and the  Chebyshev interpolation  of $T(\cdot)$.   Note
that since  the eigenvalues of (\ref{eq:Helmholtz3D})  with homogenous
Dirichlet  boundary  condition  are  real we  can  use  the  Chebyshev
interpolation technique  which target situations when  the eigenvalues
of interest lie in an  interval. Figure \ref{fig:RatChebApp} shows the
errors of each approximation versus the order of approximation $m$ for
both circular  and elliptic contour.   For simplicity, all  the errors
are evaluated  on a fine mesh  in $[-1,1]$, since arbitrary  curves in
the complex plane can be  parametrized using this interval.  From this
figure,  we can  easily see  that the  errors in  the Cauchy  rational
approximation and Chebyshev interpolation decay exponentially with the
order of the  approximation $m$, which implies that a  moderate $m$ is
usually  sufficient  to   reach  a  good  accuracy.   To  capture  the
eigenvalues  of  interest,  we  first  consider  a  circle  of  radius
$r  = 3.5$  centered at  $c  = 8.5$.   We  can then  solve the  linear
eigenvalue  problem  \nref{eq:sc} associated  with  \nref{eq:GenNLEVP}
with $m=25$ trapezoidal  quadrature nodes by performing  as many steps
of shift-and-invert  Arnoldi algorithm as  needed to extract  the $17$
eigenvalues closest to  the center $c$.  The left hand  side of Figure
\ref{fig:Cube_Dirichlet} presents  the eigenvalues computed  by Cauchy
approximation  and  those computed  by  Chebyshev  interpolation on  a
uniform mesh with $864$ triangles. Note  that, in order to make a fair
comparison between the two methods, the number of interpolation points
for Chebyshev  interpolation method is  chosen to  be the same  as the
number of  quadrature nodes  $m$.  For Chebyshev  interpolation method
the real interval  enclosing the eigenvalues of interest  is chosen as
$[5,12]$.  The  right  hand side  of  Figure  \ref{fig:Cube_Dirichlet}
illustrates  the  comparison  between  the accuracy  of  the  rational
approximation  and the  Chebyshev  interpolation. The  accuracy of  an
eigenpair  $(\lambda,  u)$  is   measured  by  the  relative  residual
$\|T(\lambda)u\|_2/\|u\|_2$.  Note  that the accuracy of  the rational
approximation  can  be  considerably  improved by  using  an  elliptic
contour      instead      of       a      circle.      A      rational
approximation~\nref{eq:TzApp}  is  then  built  using  an  elliptic
contour centered  at $c =  8.5$ with semi-major  axis $r_x =  3.5$ and
semi-minor  axis  $r_y   =  0.1$.   The  left  hand   side  of  Figure
\ref{fig:ellipse_BEM} presents the eigenvalues  computed by the Cauchy
rational   approximation  and   those   computed   by  the   Chebyshev
interpolation, and the right hand side of Figure \ref{fig:ellipse_BEM}
compares the relative residuals of the two methods.

\begin{figure}
	\centering
	\includegraphics[width=0.49\textwidth]{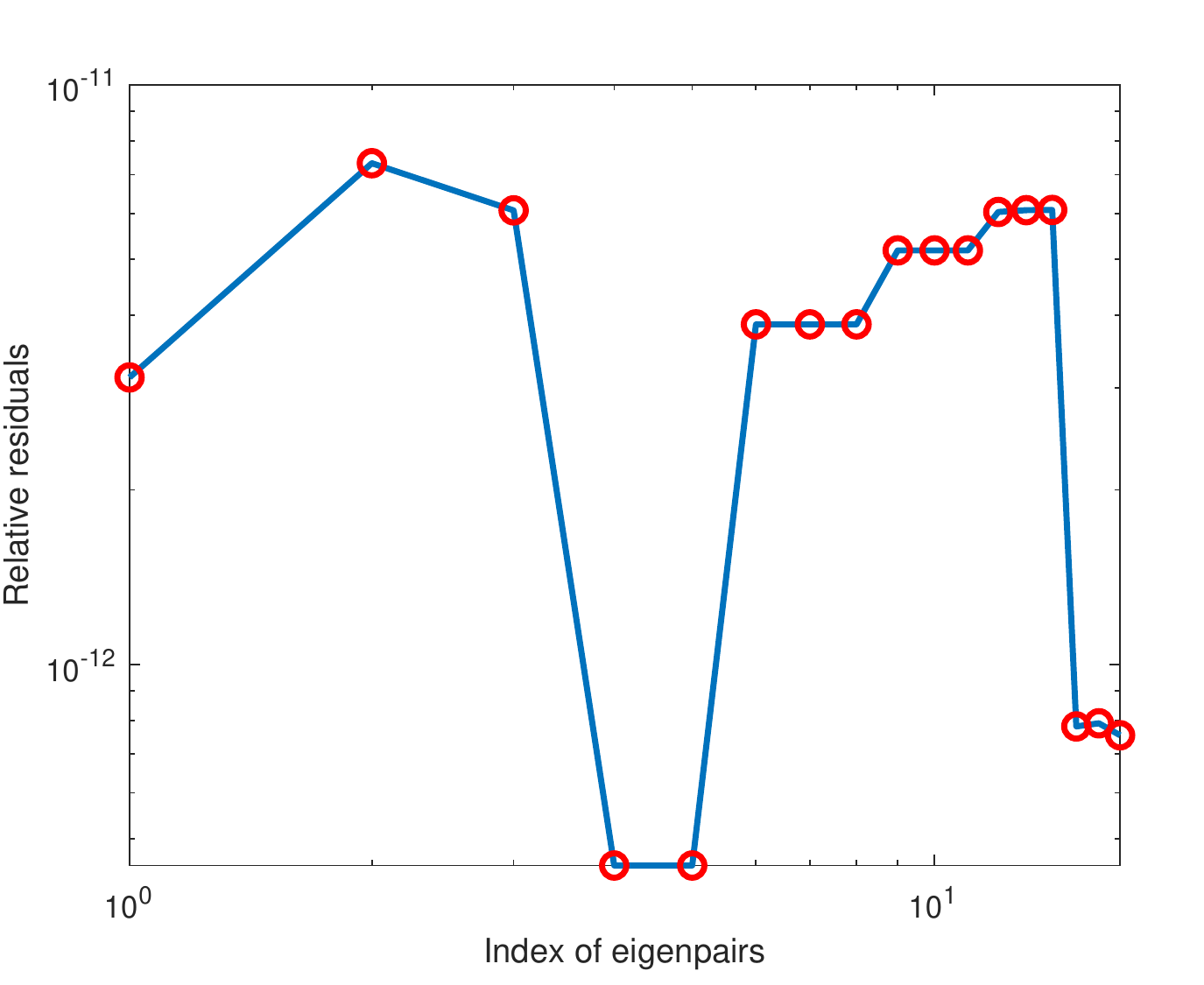}
	\includegraphics[width=0.49\textwidth]{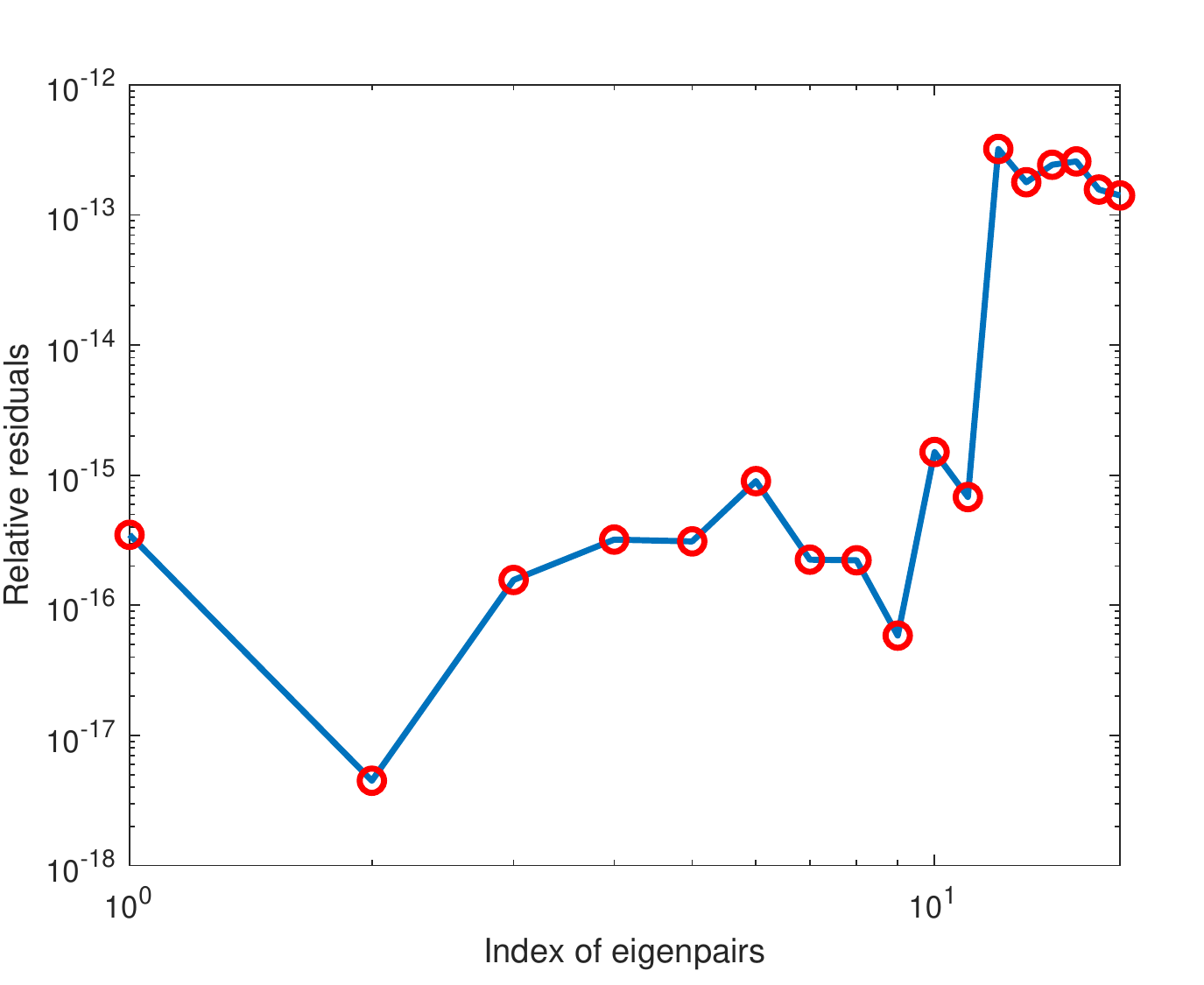}
	\caption{ Relative residuals $\|T(\lambda)u\|_2/\|u\|_2$ of the $17$ eigenvalues of
		the Laplace eigenvalue problem on the unit cube. Left: after $10$ outer iterations of the reduced approach using Cauchy approximation.
		Right: after $7$ outer iterations using Chebyshev approximation.}\label{fig:ResidualsRR}
\end{figure}

We now  repeat the same  experiment using the  Rayleigh-Ritz procedure
for  the Cauchy  and Chebyshev  approximations. To extract the  $17$
eigenvalues listed  in Table  \ref{tab:LaplaceEvals}, we start  with a
random subspace  $W$ of dimension $\nu=  20$.  We then apply  $q = 10$
steps of inverse power method to  $W$ to build a subspace of dimension
$\nu$, where each  column vector is of size $n$.  We recall that these
vectors are the  resulting top parts and bottom parts  of the iterates
of Algorithm \ref{alg:InvItCauchy} and  \ref{alg:InvItChebyshev} for the
Cauchy and Chebyshev approximations, respectively.  The resulting  subspace $U$  is then orthogonalized to
obtain  an  orthonormal  basis  $U$   that  can  be  used  to  perform
Rayleigh-Ritz projection  that leads  to a small  nonlinear eigenvalue
problem of size $\nu$. This small  problem is then solved by computing
the  eigenvalues and  eigenvectors of  the expanded  linear eigenvalue
problem \eqref{eq:EVP_red}  of size $(m+1)\nu$.  The  outer iterations
of the reduced procedure for the Cauchy approximation are stopped when
\[
\|B_0Xf_1(\Lambda)+B_2Xf_2(\Lambda)+ \ldots +B_mXf_m(\Lambda)\|_F\leq tol,
\]
and for the Chebyshev approximation when
\[
\|B_0X\tau_1(\Lambda)+B_2X\tau_2(\Lambda)+ \ldots +B_mX\tau_m(\Lambda)\|_F\leq tol,
\]
where $X$, $\Lambda$ are the extracted eigenpairs at each iteration, $\|\cdot\|_F$ denotes the Frobenius norm and $tol$ the desired tolerance for the convergence.
In our experiments, we run as many outer iterations as needed to achieve convergence with a tolerance $tol=10^{-12}$ 
for both Cauchy and Chebyshev approximations. This tolerance is achieved after $10$ outer iterations for the Cauchy approximation and after
$7$ outer iterations for Chebyshev approximation. Furthermore, Figure \ref{fig:ResidualsRR} presents
the relative residuals $\|T(\lambda)u\|_2/\|u\|_2$ for the $17$ computed eigenvalues obtained using each approximation method.
We ephasise that the Rayleigh-Ritz approach combined with Cauchy and Chebyshev approximations delivers more accurate eigenpair approximations 
than those computed by solving the linearized problem obtained directly from the Cauchy and Chebyshev approximation with projection.


\subsection*{Example 2}
As a second example, we consider the 3D Laplace eigenvalue problem \nref{eq:Helmholtz3D} on a unit sphere with homogeneous
Dirichlet boundary conditions. The analytic expressions for the eigenvalues for this geometry are well-known and given as the zeros of the spherical
Bessel function of order $\ell$. We are interested in the $6$ eigenvalues of \nref{eq:Helmholtz3D} listed in Table~\ref{tab:LaplaceSphereEvals}.

\begin{figure}
	\centering
	\includegraphics[width=0.49\textwidth]{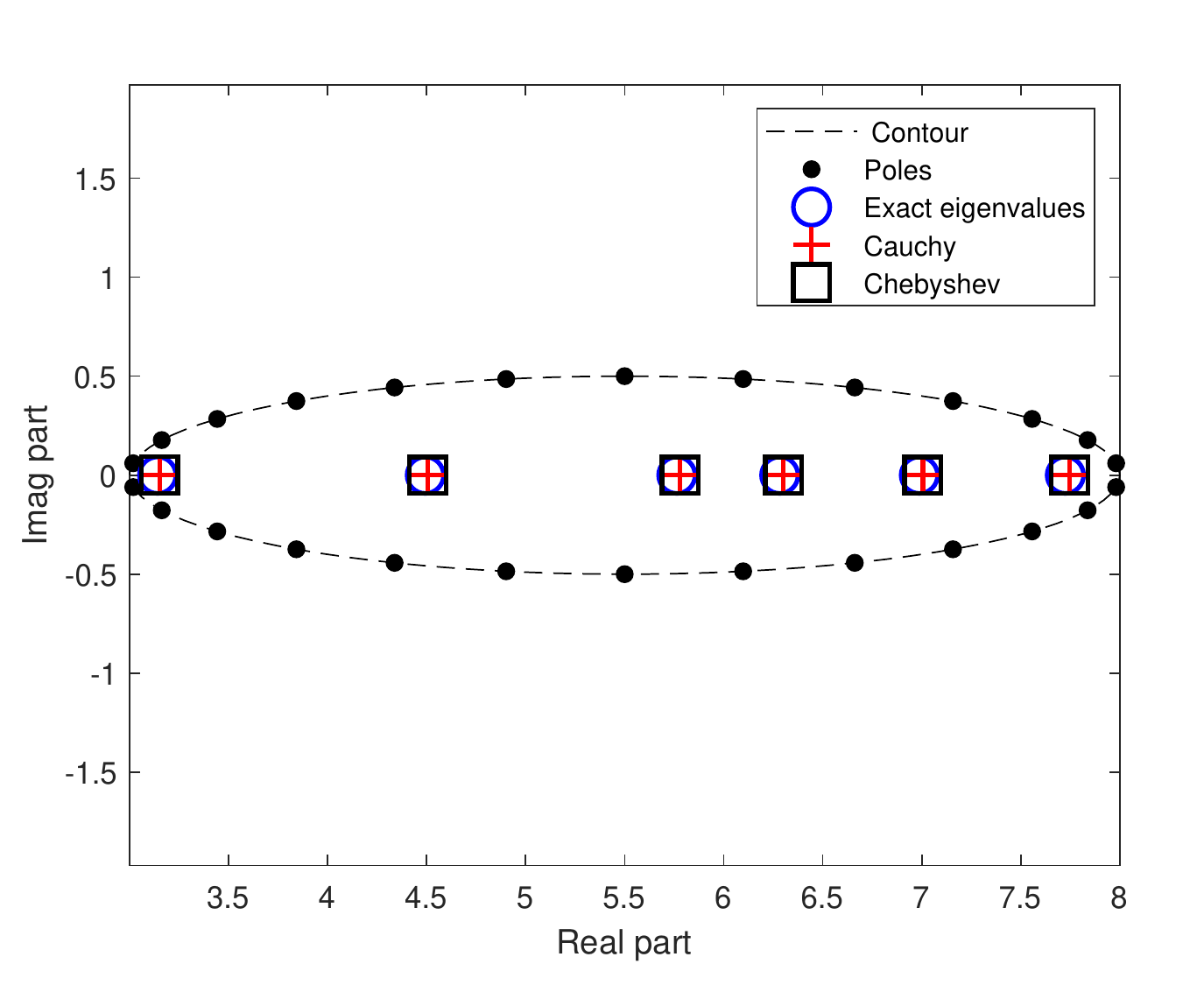}
	\includegraphics[width=0.49\textwidth]{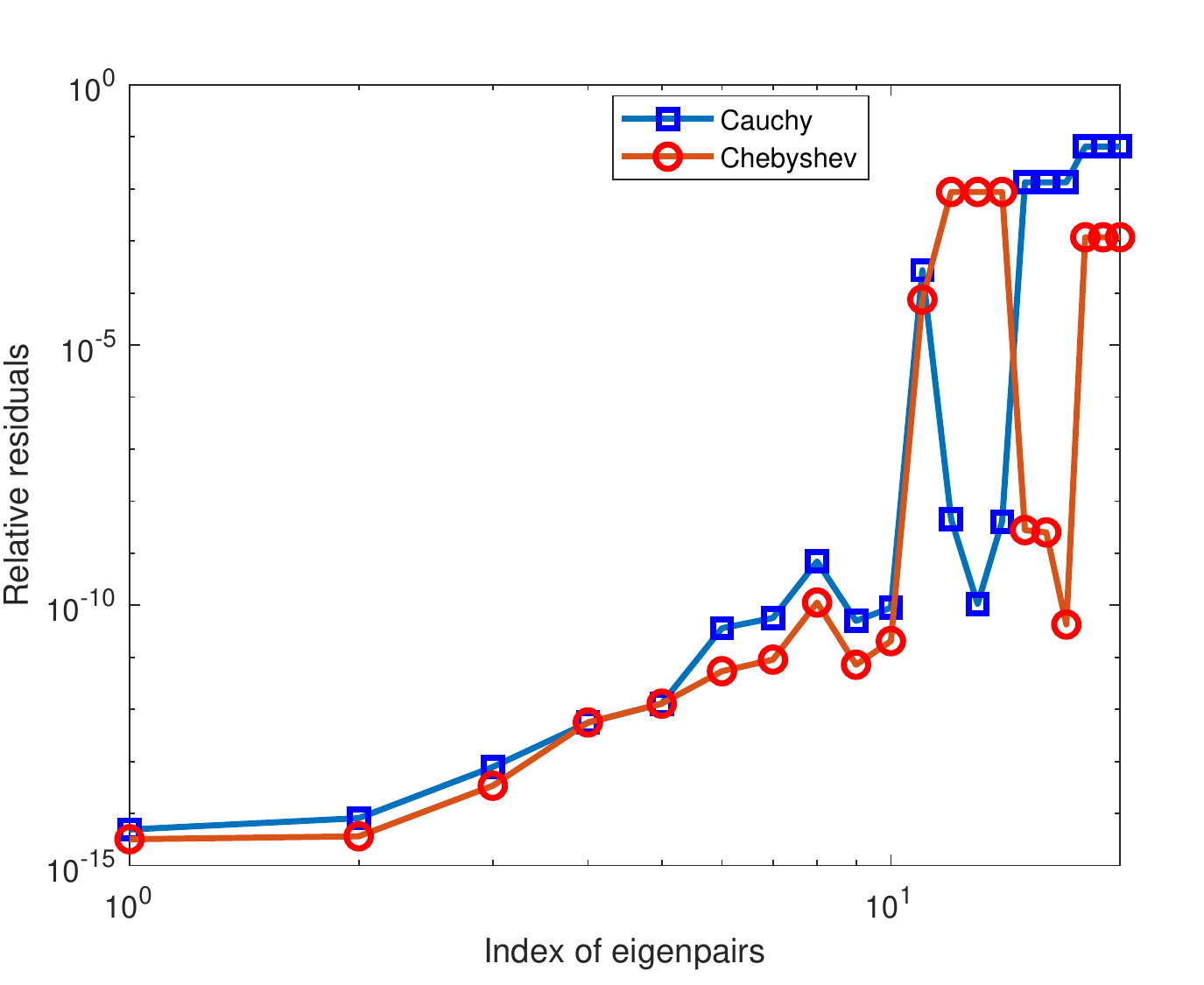}
	\caption{Left: The eigenvalues of \nref{eq:Helmholtz3D} with homogeneous Dirichlet boundary conditions inside an ellipse 
		centered at $c = 5.5$ with semi-major axis $r_x = 2.5$ and semi-minor axis $r_y = 0.2$ (circles) 
	        computed via \nref{eq:sc} (plus) and Chebyshev interpolation method inside the real interval $[3,8]$ (squares).
		Right:  The relative residuals $\|T(\lambda)u\|_2/\|u\|_2$ associated with the computed eigenpairs.}\label{fig:ellipse_BEM_sphere}
\end{figure}
\begin{figure}
	\centering
	\includegraphics[width=0.49\textwidth]{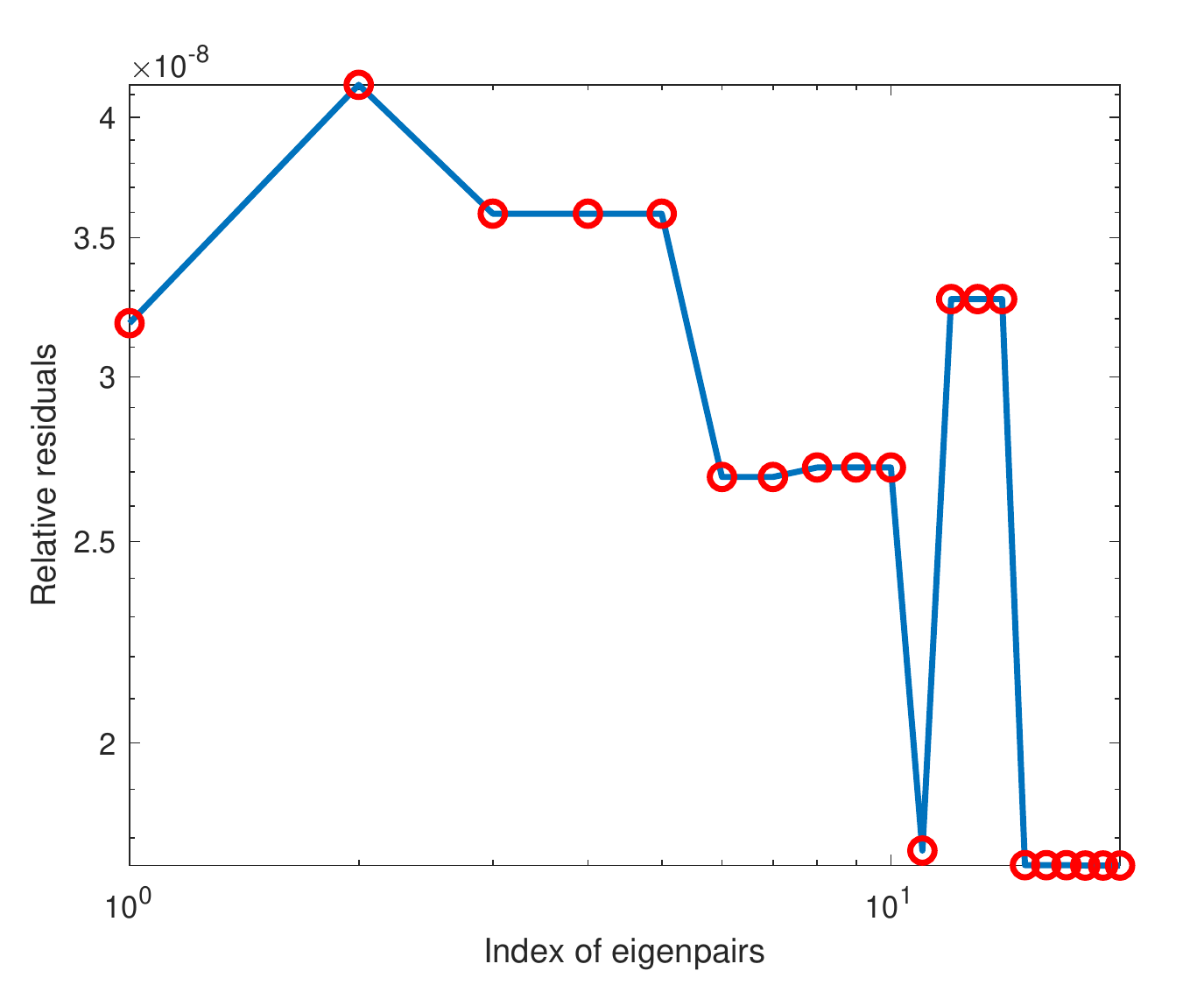}
	\includegraphics[width=0.49\textwidth]{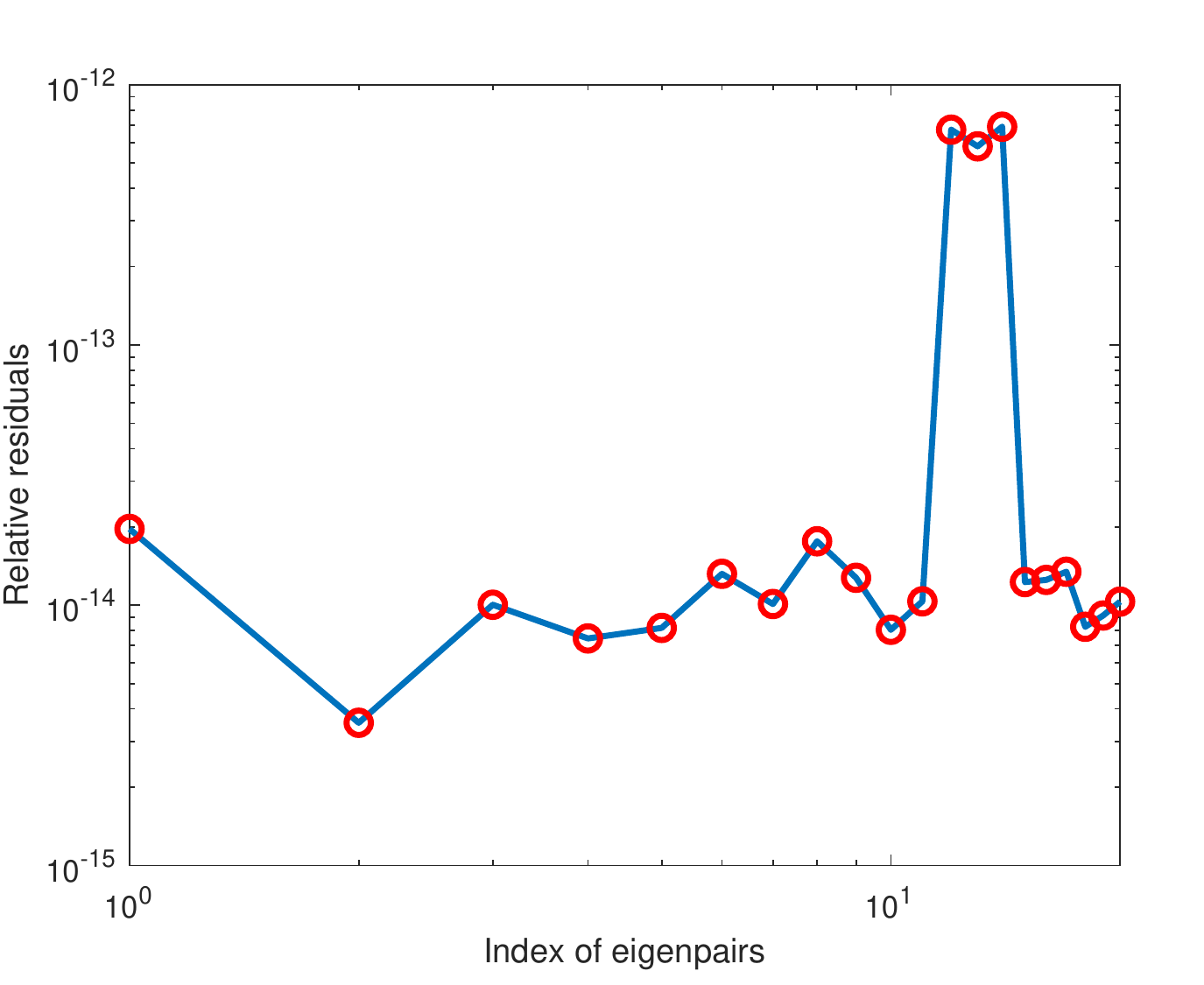}
	\caption{Relative residuals $\|T(\lambda)u\|_2/\|u\|_2$ associated with the $20$ eigenvalue approximations of
		the Laplace eigenvalue problem on the unit sphere. Left: after $17$ outer iterations of the reduced approach using Cauchy approximation. 
		Right:  after $11$ outer iterations using Chebyshev interpolation.}\label{fig:ResidualsRRsph}
\end{figure}

\begin{table}[H]
	\begin{center}
		\begin{tabular}{c c c}
			\hline
			no. & eigenvalue & multiplicity\\
			\hline
			1 & 3.1416  &1  \\
			2 & 4.4934  & 3 \\
			3 & 5.7634  & 5 \\
			4 & 6.2831 & 1 \\
			5 & 6.9879 &  7\\
			6 & 7.7252 & 3 \\
			\hline\\
		\end{tabular}
	\end{center}
	\caption{Exact eigenvalues  of the 3D Laplace eigenvalue problem on a unit sphere with homogeneous Dirichlet boundary conditions.} 
	\label{tab:LaplaceSphereEvals}
\end{table}
In  order  to compute the   eigenvalues  of  interest  using the rational
approximation technique,  we consider an elliptic  contour centered at
$c  = 5.5$  with  semi-major  axis $r_x  =  2.5$  and semi-minor  axis
$r_y = 0.1$. The left-hand side of Figure \ref{fig:ellipse_BEM_sphere}
presents the eigenvalues computed by the Cauchy rational approximation
and those  computed by  the Chebyshev interpolation  with $m=25$  on a
uniform  mesh with  $384$ triangles,  whereas the  right-hand side  of
Figure  \ref{fig:ellipse_BEM_sphere} illustrates  the accuracy  of the
two  methods.   Also in  this  example,  we  have tested  the  reduced
subspace iteration given by Algorithm \ref{alg:subsit}. To extract the
$20$ eigenvalues of interest, we  consider the Cauchy approximation on
a circle centered at $c = 5.5$  with radius $r = 2.5$. The first outer
iteration  was carried  out with  a random  subspace $W$  of dimension
$\nu= 25$ to  which $q = 10$  steps of inverse power  method, given in
Algorithm \ref{alg:subsit}, were applied. As  for Example 1, we run as
many  outer  iterations  as  needed  to  achieve  convergence  with  a
tolerance  $tol=10^{-12}$ for  both approximation  methods. Note  that
$q = 10$ steps of the inverse power  method were applied  at each outer
iteration.  The Cauchy  approximation and Chebyshev  interpolation methods
achieved  desired tolerance  after   $17$  and  $11$  outer  iterations,
respectively.   For   completeness,   we  have   also   computed   the
corresponding relative  residuals $\|T(\lambda)u\|_2/\|u\|_2$  for the
resulting eigenpairs.   These relative  residuals are shown  in Figure
\ref{fig:ResidualsRRsph} for each eigenpair index.
\begin{figure}
	\centering
	\includegraphics[width=0.49\textwidth]{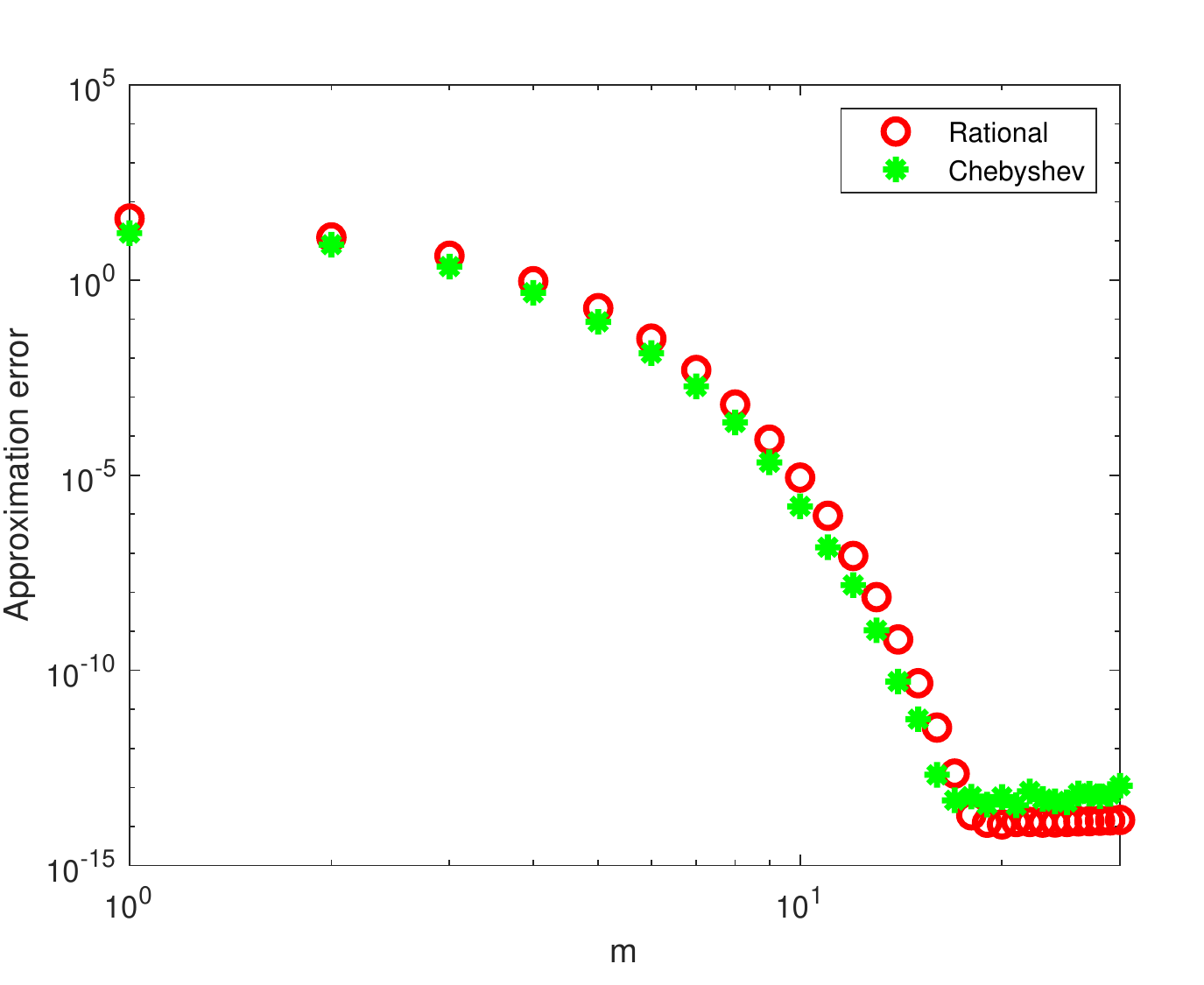}
	\caption{ Approximation errors versus the order of the approximation $m$ inside an ellipse centered at $c =0$ with semi-major axis
		$r_x = 1$ and semi-minor axis $r_y = 0.2$ for the spherical BEM problem.}\label{fig:RatChebAppSphere}
\end{figure}
\subsection*{Example 3}
\begin{figure}
	\centering
	\includegraphics[width=0.70\textwidth]{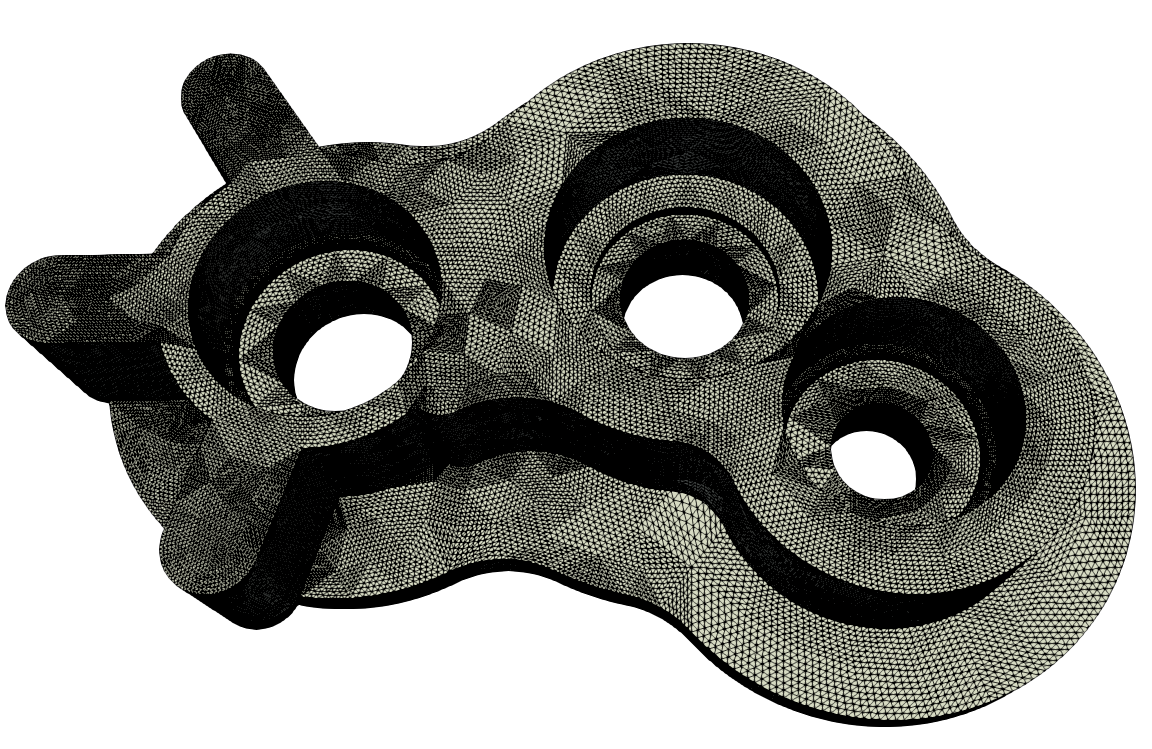}
	\caption{Geometry and BE mesh of the thermal model of a pump casing with $1 \ 728 \ 508$  DoFs.}\label{fig:pump_BEM}
\end{figure}

In  this  example,   we  illustrate  the  efficiency   of  the  Cauchy
approximation technique  applied to  the nonlinear  eigenvalue problem
resulting  from  BE  discretization   of  a real-world  problem  of
industrial  relevance. We  consider  the geometry  corresponding to  a
a pump casing model created   by  using  the   Gmsh  tool~\cite{Gmsh}. Several methods have been proposed in the literature  to comprehensively study
the acoustic behaviors of the pump casing ~\cite{pump_casing1, pump_casing2}. 
The boundary  of the  pump model  displayed in  Figure
\ref{fig:pump_BEM}  is partitioned  into $3  \ 479  \ 652$  triangles,
leading to a nonlinear eigenvalue problem with $1 \ 728 \ 508$ DoFs.
Problems of such large size add another level of difficulty to our methods, for example, we are unable to store the underlying
matrices $B_i$ in memory. To overcome this we resort to the $\mathcal{H}$-matrix based compression techniques. Specifically, we will use the {\sc Gypsilab} toolbox library {\sc openHMX}~\cite{openHMX}
in order to directly assemble $\mathcal{H}$-matrix compressed versions of matrices $B_i$. Here, we consider the boundary element discretization of problem \nref{eq:Helmholtz3D} with a rigid boundary,
i.e., $b(x)=0$. For this problem, as mentioned before, complex eigenvalues may occur. Hence, only the Cauchy approximation can be used to determine the eigenvalues
associated with this problem. Since for this example the analytic expressions of the eigenvalues are not available, the relative residuals of the computed eigenpairs
will be used to verify the accuracy of the obtained approximations. 

\begin{figure}
	\centering
	\includegraphics[width=0.49\textwidth]{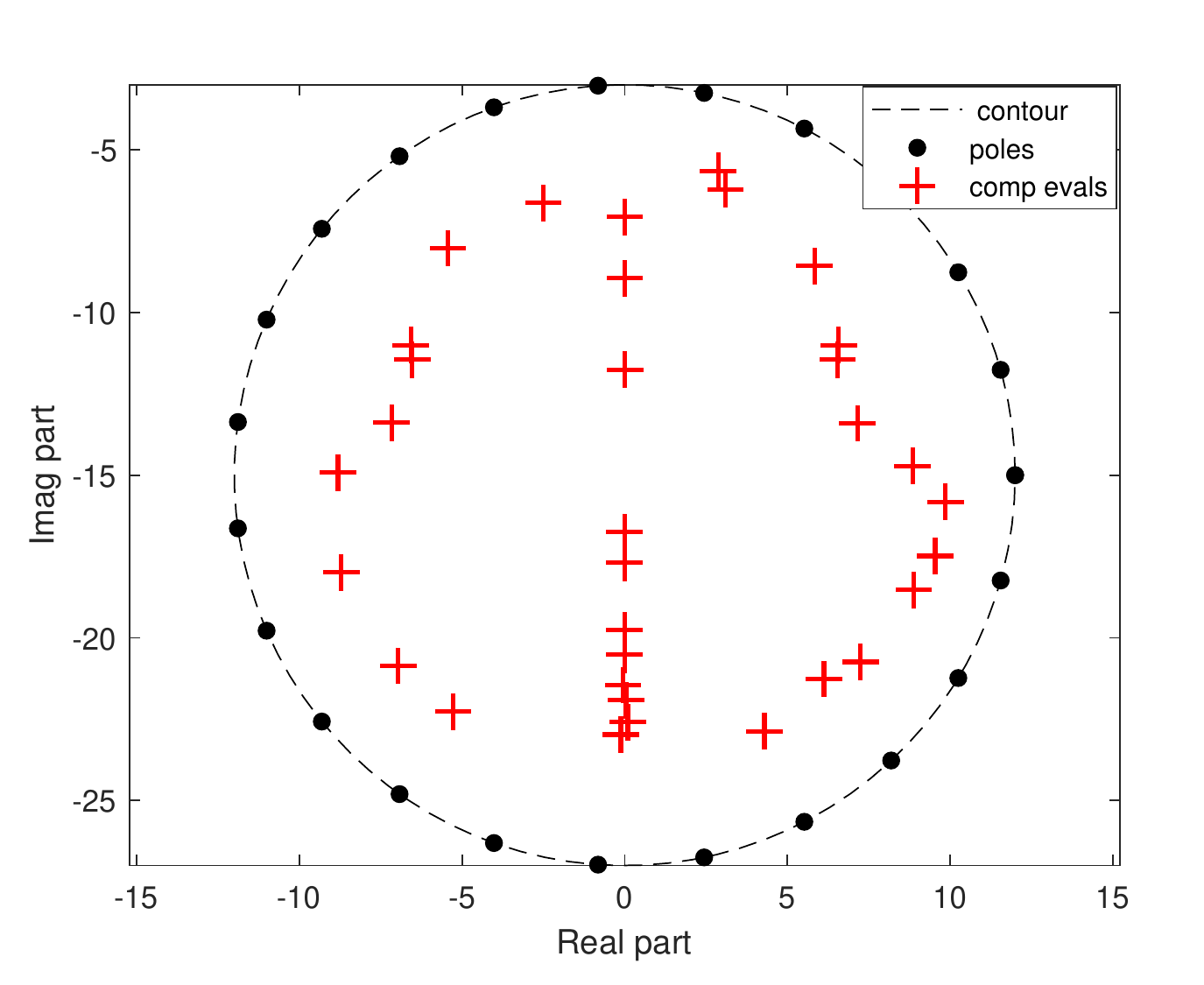}
	\includegraphics[width=0.49\textwidth]{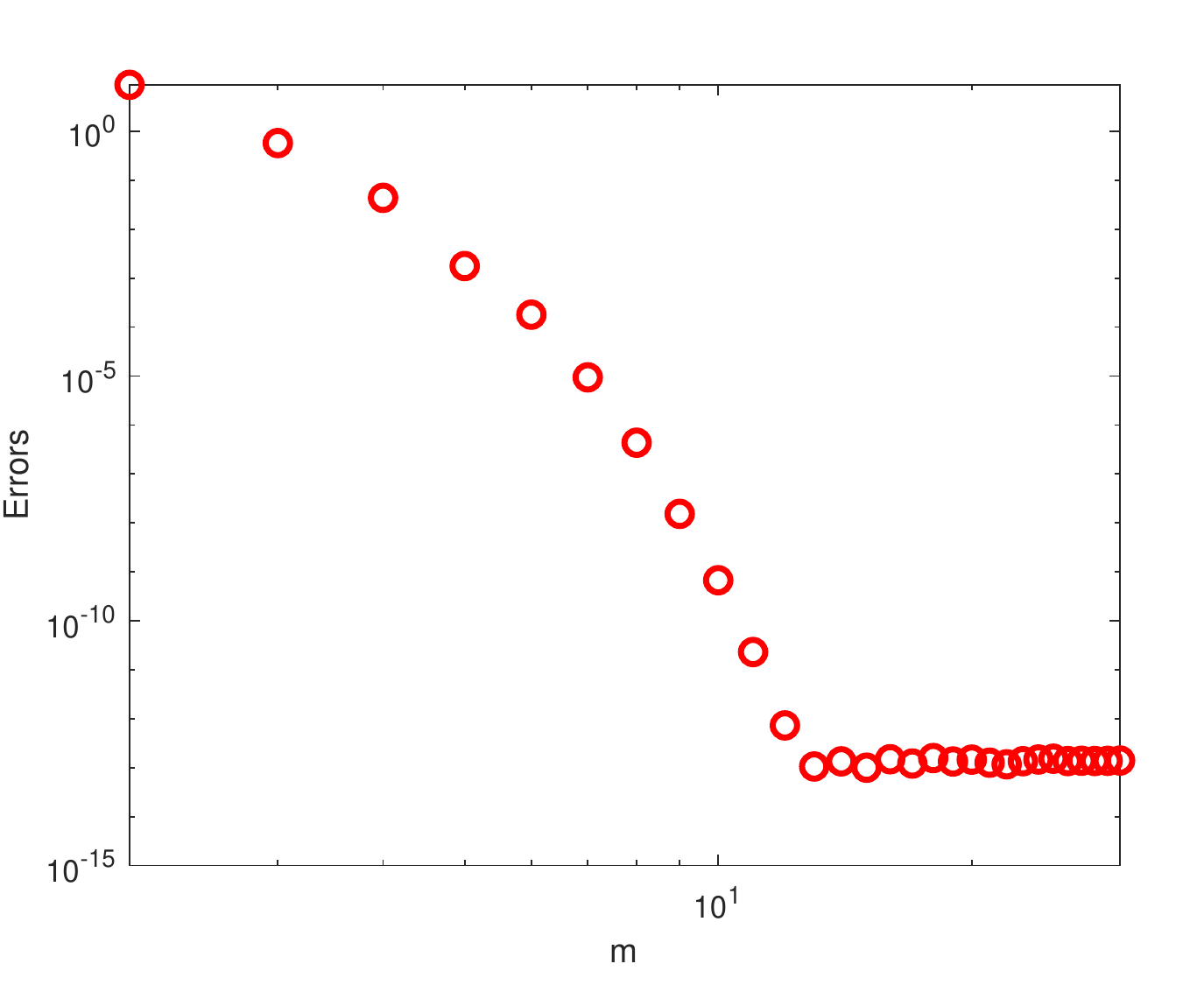}
	\caption{Left: The eigenvalues (inside a circle centered at $c = -15i$ with radius $r = 12$) of \nref{eq:Helmholtz3D} on the pump model domain with Robin boundary conditions computed via \nref{eq:sc}.
	 Right: Approximation errors versus the order $m$ of the Cauchy approximation inside a unit circle.}\label{fig:pump_robin}
\end{figure}

\begin{figure}
	\centering
	\includegraphics[width=0.49\textwidth]{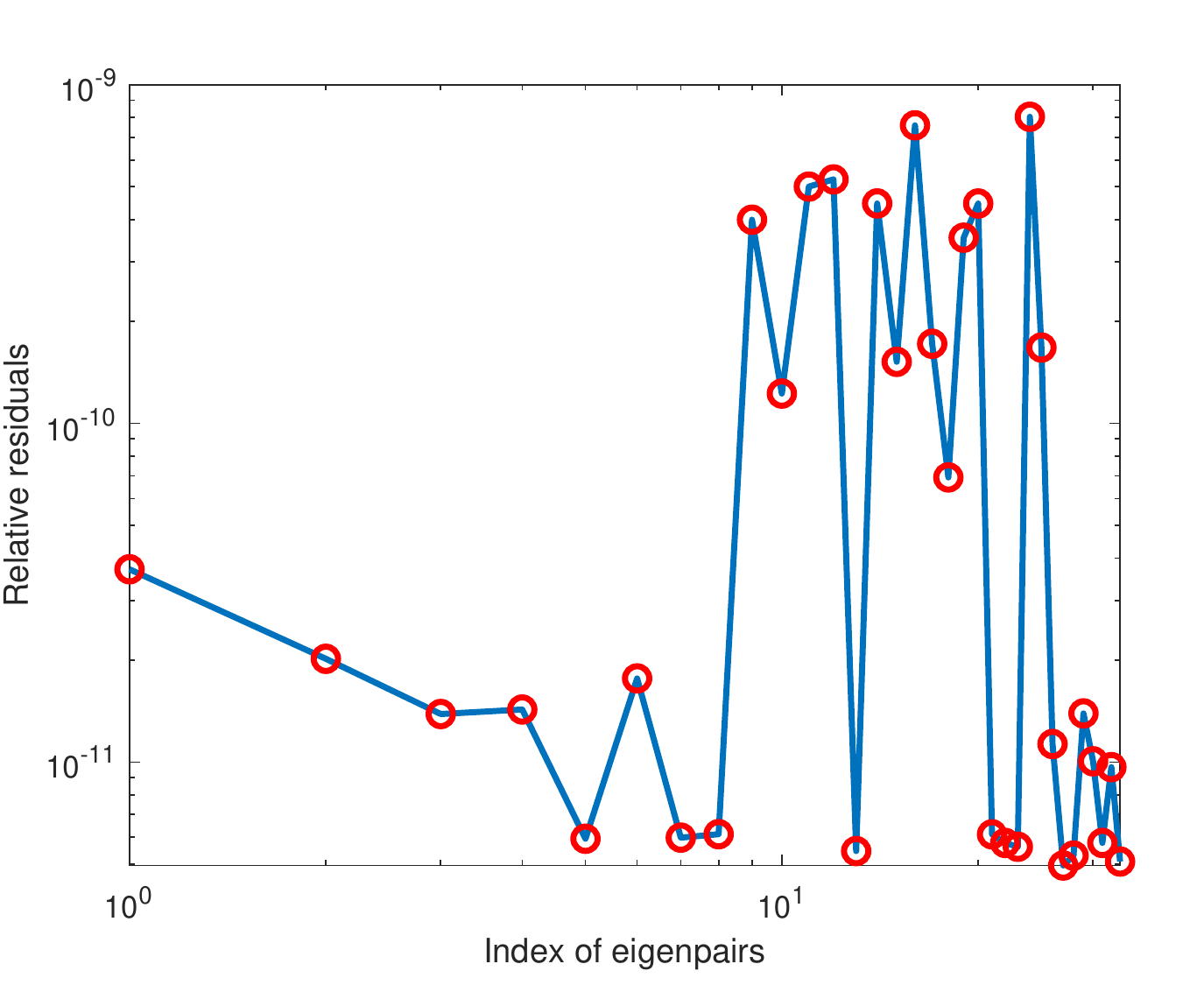}
	\caption{ Relative residuals $\|T(\lambda)u\|_2/\|u\|_2$ of the eigenvalue approximations (inside a circle centered at $c = -15i$ with radius $r = 12$) 
	of the Laplace eigenvalue problem (with Dirichlet boundary conditions) associated with the pump model displayed in Figure~\ref{fig:pump_BEM}. }\label{fig:pump_residuals}
\end{figure}

\begin{figure}[htbp]
\subfloat[Eigenfrequency: 2.87 - 5.64i Hz ]{\includegraphics[height=1.3in]{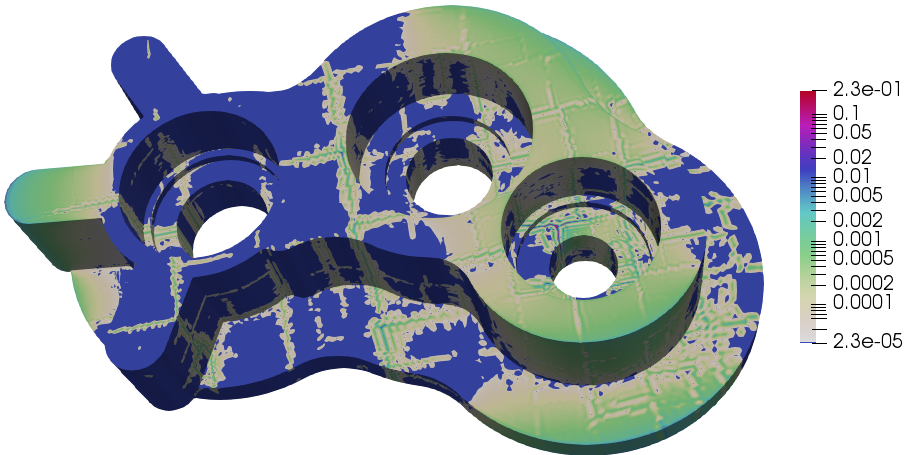}}
\subfloat[Eigenfrequency: 3.09 - 6.22i Hz]{\includegraphics[height=1.2in]{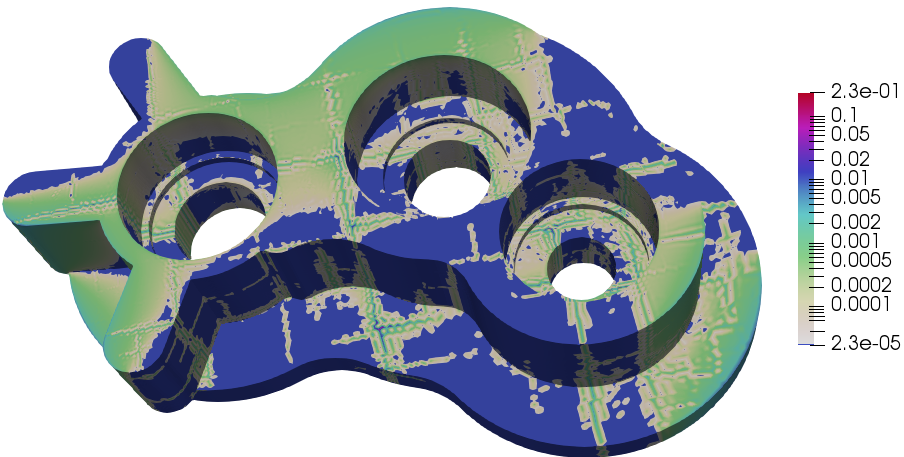}}\\
\subfloat[Eigenfrequency: 5.83 - 8.56i Hz]{\includegraphics[height=1.3in]{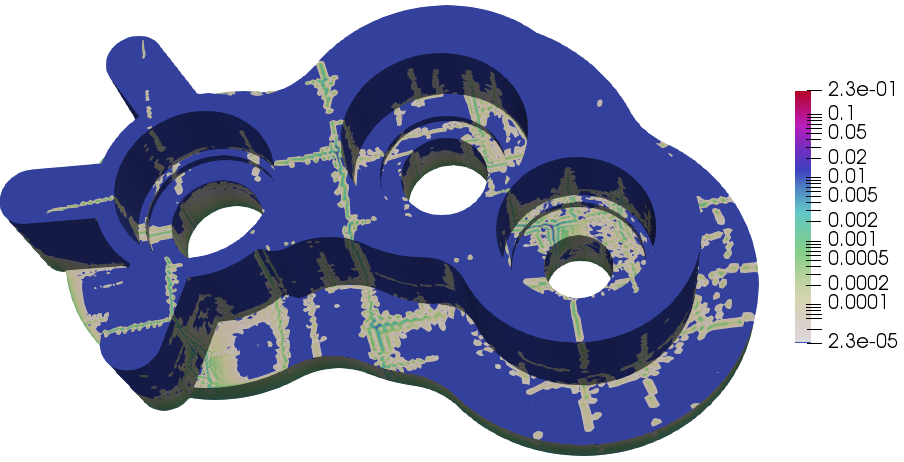}}
\subfloat[Eigenfrequency: -2.50 - 6.63i Hz]{\includegraphics[height=1.3in]{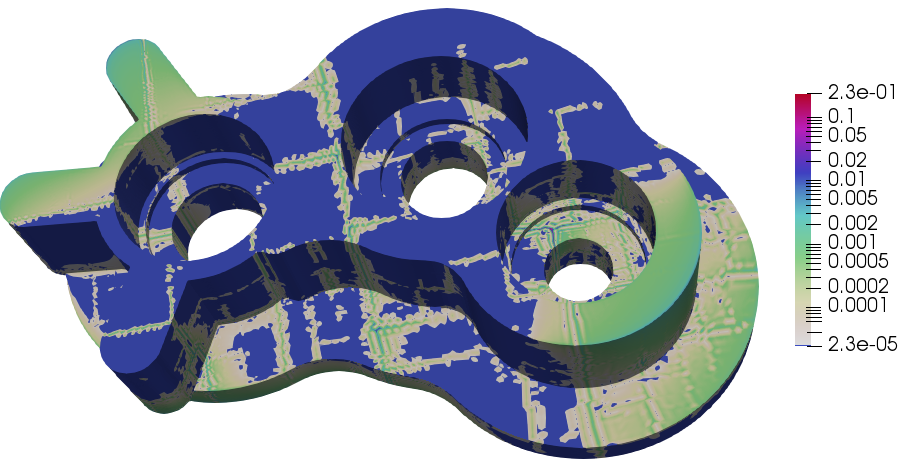}}
\caption{Eigenmodes corresponding to four different eigenvalues}
\label{fig:pump_eigenmodes}
\end{figure}

Let the  domain for the  Cauchy approximation  be given as  a circular
contour $\Omega$ centered at $c =  -15i$ with radius $r=12$. To choose
a  suitable  order  of   approximation,  we  consider  another  circle
$\Omega_1$ inside $\Omega$  with the same center  and radius $r_1=r/2$
and  then  increase $m$  until  the  resulting rational  approximation
inside $\Omega_1$  is accurate enough.  The  eigenvalue approximations
inside  $\Omega$  can be  obtained  using  a different,  much  coarser
triangular  mesh and  running as  many steps  of Arnoldi  algorithm as
needed  to accurately  solve  the expanded  linear eigenvalue  problem
\nref{eq:sc}. We recall that only  one $LU$ factorization of the Schur
complement  $\mathcal{H}$-matrix $S$  is required  in a  preprocessing
step   before   the   actual   Arnoldi  algorithm   is   invoked.   In
Figure~\ref{fig:pump_robin},  we  present   the  approximation  errors
versus  the order  of  the  Cauchy approximation  on  a  fine mesh  on
$\Omega_1$.  The right-hand side  of Figure \ref{fig:pump_robin} shows
that a high accuracy of the  rational approximation can be reached for
$m   =  24$.   We   can  therefore   solve   the  eigenvalue   problem
\nref{eq:NEVPHelmholtz}   using  $m   =  24$   trapezoidal  quadrature
nodes.   Forming  the   $24$   matrices  $B_i$   and  performing   the
matrix-vector multiplications with  $B_i$ are efficiently parallelized
on $32$ cores, where the per-core memory limit is $\approx 31$ GB. The
overall  computational time  is $7.3$  hours.  The  left-hand side  of
Figure \ref{fig:pump_robin}  shows the computed eigenvalues.  It turns
out that there are $33$  eigenvalues inside the contour $\Omega$.  The
relative  residuals  $\|T(\lambda)u\|_2/\|u\|_2$ associated  with  the
computed eigenvalues are  presented in Figure~\ref{fig:pump_residuals}
and Figure \ref{fig:pump_eigenmodes} shows  $4$ different modes of the
pump  model.

We now consider  the reduced subspace iteration approach  to solve the
same  problem  with  $10704$  triangles.  To  extract  the  same  $12$
eigenvalues   displayed    in   the    left-hand   side    of   Figure
\ref{fig:pump_robin},  we start  with a  random subspace  $W$ of  size
$\nu=20$   and  carry   out   $20$  outer   iterations  of   Algorithm
\ref{alg:subsit} with $q=10$ steps  of inverse power method (Algorithm
\ref{alg:InvItCauchy}) performed at each  single outer iteration. Note
that Algorithm  \ref{alg:subsit} is suitable for  parallelization.  In
our tests, we have exploited  parallelism for  computing
the matrices $B_i$ and  $\widehat{B}_i=U^TB_iU$ using $32$ cores. Since
Algorithm \ref{alg:InvItCauchy}  is applied to each  vector separately
to  obtain a  block of  vectors $U$,  the construction of  the approximate
subspace  at   each  outer   iteration  can   also  be   performed  in
parallel. The relative residuals reached at  the end of the $20$ outer
iterations  of  Algorithm  \ref{alg:subsit} are  displayed  in  Figure
\ref{fig:pump_residuals}.  Figure \ref{fig:pump_eigenmodes}  shows $4$
modes of the model.



\medskip
\paragraph{Acknowledgements}
The authors would like to thank Mohammed Seaid for providing the mesh data for
the test in Example~3 of the experiments and to
\texttt{Gmsh} team for making their software available.
Similarly, calculations for Example~3 in the experiments could not have been carried out without
the availability of the {\sc Gypsilab} toolbox library.
The authors benefitted from the hardware resources and support
from the Minnesota Supercomputing Institute


\bibliographystyle{siam} 

\bibliography{local}

\end{document}